\documentclass[11pt,a4paper,reqno]{amsart}
\usepackage{a4wide}
\usepackage{amsfonts,amsmath,amssymb,amsthm,enumerate,color,hyperref,bbm,tabularx,ifthen,twoopt,mathabx}
\usepackage{graphicx}
\usepackage[latin1]{inputenc}

\newcounter{hypA}

\def\PE{\mathbb{E}}
\def\PP{\mathbb{P}}
\def\Var{\mathop{\rm Var}\nolimits}
\def\Tr{\mathop{\rm Tr}\nolimits}
\def\rmd{\mathrm{d}}
\def \1{\mathbbm{1}}
\def\Y{\mathbf{Y}}
\def\betabold{\boldsymbol{\beta}}
\def\X{\mathbf{X}}
\def\u{\mathbf{u}}
\def\Z{\mathbf{Z}}
\def\R{\mathbf{R}}
\def\W{\mathbf{W}}
\def\D{\mathbf{D}}
\def\B{\mathbf{B}}
\def\M{\mathbf{M}}
\def\H{\mathbf{H}}
\def\U{\mathbf{U}}
\def\T{\mathbf{T}}
\def\J{\mathbf{J}}
\def\e{\mathbf{e}}
\def\v{\mathbf{v}}
\def\t{\mathbf{t}}
\def\V{\mathbf{V}}
\def\L{\mathbf{L}}
\def\w{\mathbf{w}}
\def\eps{\boldsymbol{\varepsilon}}
\def\Id{\textrm{Id}}
\def\rset{\mathbb{R}}

\newtheorem{theorem}{Theorem}
\newtheorem*{theorem2}{Theorem}

\newtheorem{lemma}{Lemma}

\theoremstyle{remark}

\newtheorem{assumption}{Assumption}

\begin{document}

\title{Heritability estimation in high dimensional linear mixed models}

\author{A. Bonnet}
\author{E. Gassiat}
\author{C. L\'evy-Leduc}
\address{Laboratoire de Math\'ematiques de l'Universit\'e Paris-Sud}
\email{elisabeth.gassiat@math.u-psud.fr}
\address{AgroParisTech/UMR INRA MIA 518}
\email{celine.levy-leduc@agroparistech.fr, anna.bonnet@agroparistech.fr}

\maketitle

\begin{abstract}
Motivated by applications in genetic fields, we propose to estimate the heritability in high dimensional sparse linear mixed models. 
The heritability determines how the variance is shared between the different random components of a linear mixed model.
The main novelty of our approach is to consider that the random effects can be sparse, that is may contain null components, 
but we do not know neither their proportion nor their positions. The estimator that we consider is strongly inspired
by the one proposed by \cite{pirinen:donnelly:spencer:2013}, and is based on a maximum likelihood approach. 
We also study the theoretical properties of our estimator, namely we establish that our estimator of the heritability
is $\sqrt{n}$-consistent when both the number of observations $n$ and the number of random effects $N$ tend to infinity
under mild assumptions. We also prove that our estimator of the heritability satisfies a central limit theorem 
which gives as a byproduct a confidence interval for the heritability. 
Some Monte-Carlo experiments are also conducted in order to show the finite sample performances of our estimator.
\end{abstract}


\section{Introduction}
\label{sec:intro}

Linear mixed models (LMM) have been widely used in various fields such as agriculture, biology, medicine and genetics.
More precisely, in quantitative genetics, LMM have been used for estimating the heritability of traits and breeding values
as explained for instance by \cite{lynch:1998}. In Genome Wide Association Studies (GWAS), which is 
the application field that inspired our work, \cite{yang:lee:goddard:visscher:2011} suggested the use of 
linear mixed models to measure genotypes at a large number of single nucleotide polymorphisms (SNP) 
-typically 300,000 to 500,000- in large sample of individuals -typically, 1000- in order to identify genetic 
variants that explain phenotypes variations.


The model that we shall study in this paper is a LMM defined as follows:
\begin{equation}\label{eq:modele}
\Y=\X\betabold+\Z\u+\e\;, 
\end{equation}
where $\Y=(Y_1,\dots,Y_n)'$ is the vector of observations, $\X$ is a $n\times p$ matrix of predictors, $\betabold$
is a $p\times 1$ vector containing the unknown linear effects of the predictors, $\u$ and $\e$ correspond to the random
effects. Moreover, in (\ref{eq:modele}), $\Z=(Z_{i,j})$ is a $n\times N$
matrix such that the $Z_{i,j}$ are normalized random variables in the following sense: they are defined from 
a matrix $\W=(W_{i,j})_{1\leq i\leq n,\,1\leq j\leq N}$ by
\begin{equation}\label{eq:normalization_1}
Z_{i,j}=\frac{W_{i,j}-\overline{W}_{j}}{s_{j}},\;i=1,\dots,n,\;j=1,\dots,N\;,
\end{equation}
where
\begin{equation}\label{eq:normalization_2}
\overline{W}_{j}=\frac{1}{n}\sum_{i=1}^{n}W_{i,j},\;s_{j}^{2}=\frac{1}{n}\sum_{i=1}^{n}(W_{i,j}-\overline{W}_{j})^{2},
\;j=1,\ldots,N\;.
\end{equation}
In (\ref{eq:normalization_1}) and (\ref{eq:normalization_2}) the $W_{i,j}$'s are such that
for each $j$ in $\{1,\dots,N\}$ the $(W_{i,j})_{1\leq i\leq n}$ are independent and identically distributed 
random variables and such that the columns of $\W$ are independent. 
With this definition the columns of $\Z$ are empirically centered and normalized.
Finally, we shall assume that the random effects can be sparse, that is only a proportion $q$ of the components of $u$ are non zero:
\begin{equation}\label{eq:distrib_u_e_comp_nulles}
u_i\stackrel{i.i.d.}{\sim}(1-q)\delta_0+q\mathcal{N}(0,{\sigma_u^\star}^2)\;,\textrm{for all } 1\leq i\leq N \textrm{ and }
\e\sim\mathcal{N}\left(0,{\sigma_e^\star}^2\Id_{\rset^n}\right)\;,
\end{equation}
where $\Id_{\rset^n}$ denotes the $n\times n$ identity matrix, $q$ is in $(0,1]$, and $\delta_0$ is the point mass at $0$. Notice that the case $q=1$ corresponds to the usual non sparse model.
In Model (\ref{eq:modele}) with (\ref{eq:normalization_1}), (\ref{eq:normalization_2}), (\ref{eq:distrib_u_e_comp_nulles}), 
one can observe that
$$
\Var(\Y|\Z)=Nq{\sigma_u^\star}^2 \R +{\sigma_e^\star}^2 \Id_{\rset^n}\;,\textrm{ where }\R=\frac{\Z\Z'}{N}\;.
$$
Inspired by \cite{pirinen:donnelly:spencer:2013}, Model (\ref{eq:modele}) can be rewritten by using the following parameters:
\begin{equation}\label{eq:param_comp_nulles}
{\sigma^\star}^2=Nq{\sigma_u^\star}^2+{\sigma_e^\star}^2 \textrm{ and } \eta^\star
=\frac{Nq{\sigma_u^\star}^2}{Nq{\sigma_u^\star}^2+{\sigma_e^\star}^2}\;.
\end{equation}
Thus,
$$
\Var(\Y|\Z)=\eta^\star{\sigma^\star}^2 \R+(1-\eta^\star){\sigma^\star}^2 \Id_{\rset^n}\;.
$$
The parameter $\eta^\star$ which belongs to $[0,1]$ is commonly called the heritability in the case where $q=1$, 
see for instance \cite{Yang:2010}, and determines how 
the variance is shared between $\u$ and $\e$ when all the components of $\u$ are non zero. We propose in
(\ref{eq:param_comp_nulles}) to extend this definition to the case where $\u$ may contain null components and 
$q$ is in $(0,1]$.
Our main goal in this paper is to propose an estimator for the heritability in this possibly sparse framework
 and to establish its theoretical
properties in the non standard theoretical context where $n$ and $N$ tend to infinity.

 The use of the linear mixed models to estimate heritability has been proposed by \cite{yang:lee:goddard:visscher:2011} as an alternative to the regression models usually used in the GWAS. The goal is to consider the joint effect of all SNPs on a phenotype, and the heritability corresponds to the proportion of phenotypic variance
explained by all SNPs.
More precisely, the $i$th component $u_i$ of $\u$ corresponds to the effect of the $i$th SNP on the phenotype
and $\e$ corresponds to the environmental effect. 
In this application, the matrix $\W$ contains all the genetic information about all the individuals in the study. More precisely, for each $j$, the $(W_{i,j})_{1\leq i\leq n}$ are i.i.d binomial random variables with parameters 2 and $p_j$. $W_{i,j}=0$ (resp. 1, resp. 2) if the genotype of the $i$th individual at locus $j$ is $qq$
(resp. Qq, resp. QQ) where $p_j$ is the frequency of Q allele at locus $j$. In the GWAS framework, $\Z$ is thus a matrix having 
a number of rows equal to the number of individuals in the experiment that is $n\approx 1000$ and a number of columns
equal to the number of SNPs taken into account in the experiment, namely $N\approx 500,000$. This application motivated the framework
that we have chosen where $n$ and $N$ tend to infinity and where $Z$ is a random matrix.

The major difference between the framework of  \cite{yang:lee:goddard:visscher:2011} and ours is that they consider that the random effects are 
Gaussian while we consider a mixture model between a point mass at 0 and a Gaussian distribution.
 With this modeling, we assume that all SNPs are not necessarily causal, that is that all SNPs do not explain a given phenotype.
The parameter $q$ defined in (\ref{eq:distrib_u_e_comp_nulles}) actually corresponds to the proportion of non null components in $\u$ that is to the proportion of causal SNPs. Then, the heritability defined by $\eta^\star$ in (\ref{eq:param_comp_nulles}) corresponds 
to the proportion of phenotypic variance explained by the causal variants.

Several approaches can be used for estimating the heritability in the case where $q=1$ but 
to the best of our knowledge, no theoretical results concerning the estimation of the heritability or the estimation
of $\sigma_u^\star$, $\sigma_e^\star$ have been established in the framework where both $n$ and $N$
tend to infinity. This is one of the contributions of our paper. Among these approaches, we can quote
the REML (REstricted Maximum Likelihood) approach, originally proposed by \cite{patterson:thompson:1971} and
described for instance in \cite{searle:casella:mcculloch:1992}, which consists in estimating $\sigma_u^\star$ and $\sigma_e^\star$
for estimating $\eta^\star$. However, this type of approach is based on iterative procedures that require expensive
matrix operations. Hence, several approximations have been proposed such as the \textsf{AI} algorithm 
(\cite{gilmour:thompson:cullis:1995}) which is used for instance in the software GCTA (Genome-wide Complex Trait Analysis)
described in \cite{yang:lee:goddard:visscher:2011}. Other approximations have also been proposed in the EMMA algorithm 
(\cite{Kang:2008}). For further details on the different approximations that could be used we refer the reader to 
\cite{pirinen:donnelly:spencer:2013}. The latter paper proposes another methodology for estimating the heritability
which consists in rewriting Model (\ref{eq:modele}) with the parameters (\ref{eq:param_comp_nulles}) and in using an
eigenvalue decomposition of the matrix $\R$. Further details on their methodology are given hereafter.
According to the numerical experiments conducted in \cite{pirinen:donnelly:spencer:2013}
their approach has the lowest computational burden among the available algorithms.

In the case of sparse high dimensional framework, most of the papers studied the case of
linear models. Among them, we can quote: \cite{meinshausen:buhlmann:2010} and \cite{blanchard:2014}.
The high dimensional linear mixed model where $\u$ is sparse, that is the case where $q<1$,
which is the most realistic case for the applications that we have in view, has 
received little attention. It has
been studied
according two directions: detection and estimation. Concerning the detection field in this framework,
we are only aware of the work of \cite{castro:candes:plan:2011} in which a testing procedure 
for detecting a sparse vector in high dimensional linear sparse regression model is also proposed
and compared with the one proposed by \cite{ingster2010}. As for the procedures dedicated to the heritability 
estimation, there exist, to the best of our knowledge, only three approaches:
the approach of \cite{Yang:2010} who propose to approximate the genetic correlation 
between every pairs of individuals across the set of causal SNPs by the genetic correlation across the set of all SNPs,
the approach of \cite{golan:rosset:2011} who propose a methodology based on 
MCEM (Monte-Carlo expectation-maximization) developed by \cite{Wei:1990} and 
the Bayesian approaches of \cite{guan:stephens:2011}
and \cite{zhou:carbonetto:stephens:2013}.
However, as far as the estimation issue in
the high dimensional linear mixed model is concerned,
the authors of these papers did not establish the theoretical properties of their estimators 
in the framework where both $n$ and $N$ tend to infinity.

In this paper, we prove that
using a strategy close to the one proposed by \cite{pirinen:donnelly:spencer:2013},
which has been devised in the case $q=1$, provides consistent estimators even in the case where $q<1$.
Moreover, we prove that this 
estimator is $\sqrt{n}$-consistent in the following asymptotic framework:
$n\to\infty$ and $N\to\infty$ such as $n/N\to a >0$ and satisfies under mild assumptions a central limit theorem
in both cases $q=1$ and $q<1$. 
It has to be noticed that the classical results that exist in linear mixed models
are established only in the case where $q=1$, $n$ tends to infinity and $N$ is constant.
Note that in our asymptotic framework where $\eta^\star$ is a constant and
$N$ tends to infinity, $\sigma_u^\star$ tends to 0 as $N$ tends to infinity.

In the sequel, up to considering the projection of $\Y$ onto the orthogonal of the image of $\X$ and for notational simplicity, 
we shall focus on the following model
\begin{equation}\label{eq:modele_comp_nulles}
\Y=\Z\u+\e\;, 
\end{equation}
where the assumptions on $\u$ and $\e$ are given in (\ref{eq:distrib_u_e_comp_nulles}).

In the case where $q=1$, observe that 
$$
\Y|\Z\sim\mathcal{N}\left(0,\eta^\star{\sigma^\star}^2 \R+(1-\eta^\star){\sigma^\star}^2 \Id_{\rset^n}\right)\;.
$$
Let $\U$ as the orthogonal matrix ($\U'\U=\U\U'=\Id_{\rset^n}$) such that
$\U\R\U'=\textrm{diag}(\lambda_1,\dots,\lambda_n)$ is a diagonal matrix having its diagonal entries
equal to $\lambda_1,\dots,\lambda_n$.
Hence, in the case where $q=1$ and conditionally to $\Z$,
$\widetilde{\Y}=\U'\Y$ is a zero-mean Gaussian vector and covariance matrix equal to 
$\textrm{diag}(\eta^\star{\sigma^\star}^2\lambda_1+(1-\eta^\star){\sigma^\star}^2
,\dots,\eta^\star{\sigma^\star}^2\lambda_n+(1-\eta^\star){\sigma^\star}^2)$,
where the $\lambda_i$'s are the eigenvalues of $\R$.

The method proposed by \cite{pirinen:donnelly:spencer:2013} consists in computing the log-likelihood 
$$ L_n(\sigma^2,\eta)= -\frac{n}{2}\log(\sigma^2) -\frac{1}{2} \sum_{i=1}^{n} \log(\eta(\lambda_i-1)+1) -\frac{1}{2\sigma^2} \sum_{i=1}^{n}\frac{\widetilde{Y_i}^2}{\eta(\lambda_i-1)+1} -\frac{n}{2}\log(2\pi) $$
and to maximize this function of two variables by iterative optimization techniques.
Since in our case we are only interested in estimating $\eta^\star$, we plugged an estimator of $\sigma^{\star 2}$ that is 
$$\hat{\sigma}^2=\frac{1}{n} \sum_{i=1}^{n} \frac{\widetilde{Y_i}^2}{\eta(\lambda_i-1)+1}$$ in $L_n$.
Thus, in the case $q=1$, the maximum likelihood strategy would lead to estimate $\eta^\star$,  assumed to be in $[0, 1-\delta]$ with $\delta > 0$, 
by $\hat{\eta}$ defined as a
maximizer of 
\begin{equation}\label{eq:Ln}
L_n(\eta)= -\log \left(\frac{1}{n} \sum_{i=1}^{n} \frac{\widetilde{Y}_i^{2}}{\eta(\lambda_i-1)+1}\right) 
-\frac{1}{n}\sum_{i=1}^{n} \log\left(\eta(\lambda_i-1)+1\right)\;,
\end{equation}
where the $\widetilde{Y}_i$'s are the components of the vector $\widetilde{\Y}=\U'\Y$.

We shall establish in Theorem \ref{th:modele_general}, which is proved in Section \ref{sec:suppl}, 
that this strategy produces $\sqrt{n}$-consistent estimators
of $\eta^\star$ in both cases: $q=1$ and $q<1$ and also that this estimator satisfies a central limit theorem
which provides as a by-product confidence intervals for $\eta^\star$.

The paper is organized as follows. Section \ref{sec:theory} is dedicated to the theoretical properties of our estimator. 
The numerical results are presented in Section \ref{sec:numeric}. They have been obtained thanks to the R package
\textsf{HiLMM} that we have developed and which is available from the Comprehensive R Archive
Network (CRAN). In Section \ref{sec:discussion}, we provide some additional comments
on our work as well as some prospects such as the estimation of the proportion $q$ of non null components in the random effects. 
Finally, the proofs are given in Section \ref{sec:proofs}.


\section{Theoretical results}\label{sec:theory}

Observe that another way of writing Model (\ref{eq:modele_comp_nulles}) with the parameters defined in (\ref{eq:param_comp_nulles})
consists in writing
\begin{equation}\label{eq:modele_final}
\Y=\frac{1}{\sqrt{N}}\Z\t+\sigma^\star\sqrt{1-\eta^\star}\eps\;,
\end{equation}
where $\eps$ is a $n\times 1$ 
Gaussian vector having a covariance matrix equal to identity and $\t=(t_1,\dots,t_N)'$ is a random vector such 
that 
$$
t_i=\frac{\sigma^\star\sqrt{\eta^\star}}{\sqrt{q}}w_i\pi_i\;,
$$
where the $w_i$'s and the $\pi_i$'s are independent, $\w=(w_1,\dots,w_N)'$ is a Gaussian vector with a covariance matrix
equal to identity and the $\pi_i$'s are i.i.d Bernoulli random variables such that $\PP(\pi_1=1)=q$.

The estimator $\hat{\eta}$ is defined as a maximizer of $L_{n}(\eta)$ for $\eta\in [0,1-\delta]$ for some small $\delta>0$, $L_{n}$ being given in (\ref{eq:Ln}).  
We shall study the asymptotic properties of $\hat{\eta}$ as $n$ and $N$ tend to infinity in a comparable way, that is when $n/N\to a>0$. To understand the asymptotic behavior of $\hat{\eta}$, we shall first prove its consistency, then use a Taylor expansion of the derivative of $L_{n}$ around  $\hat{\eta}$ in the usual way. The computations as can be seen in  (\ref{eq:Ln}) involve empirical means of functions of the eigenvalues $\lambda_{i}$ of $\R=\frac{\Z\Z'}{N}$. 
Using Theorem 1.1 of \cite{bai:zhou:2008}, we shall prove the almost sure convergence of such empirical quantities under a weak assumption denoted by 
Assumption \ref{assum:A1} as follows.

\begin{assumption}\label{assum:A1}
Let $\Z$ and $\W$ be two matrices defined by (\ref{eq:normalization_1}) and (\ref{eq:normalization_2}). Recall that for each $j$ in $\{1,\dots,N\}$ the $(W_{i,j})_{1\leq i\leq n}$ are independent and identically distributed 
random variables and such that the columns of $\W$ are independent (but not necessarily identically distributed).
Assume that the entries $W_{i,j}$ of $\W$ are uniformly bounded, and have variance uniformly lower bounded, that is: there exist $W_{M}<\infty$ and $\kappa>0$ such that
$0\leq W_{i,j}\leq W_M$ and $\sigma_j^2=\Var(W_{i,j})\geq\kappa $, for all $j$.
\end{assumption}

The following lemma ensures that the result of \cite{marchenko:pastur:1967} which gives the empirical
spectral distribution of sample covariance matrices $\Z\Z'/N$ holds even when the entries
$Z_{i,j}$ of the matrix $\Z$ are not i.i.d. random variables but when $\Z$ is obtained by empirical standardization
of a matrix $\W$ satisfying Assumption 1.

\begin{lemma}\label{lem:MP_noniid}
Under Assumption \ref{assum:A1},  as $n,N\to\infty$ such that $n/N\to a>0$, the empirical spectral distribution of $R_N=\Z\Z'/N$:
$F^{R_N}(x)=n^{-1}\sum_{k=1}^n \1_{\{\lambda_k\leq x\}}$ tends almost surely to the Marchenko-Pastur distribution defined as the distribution function of $\mu_{a}$ where, for any measurable set $A$, 
$$
\mu_{a}(A)=\left\{\begin{array}{ll}\left(1-\frac{1}{a}\right)\1_{0\in A} + \nu_{a}(A)&\textrm{ if $a>1$}\\
\nu_{a}(A)&\text{ if $a\leq 1$}
\end{array}
\right.
$$
with
\begin{equation}\label{eq:marchenko}
\rmd\nu_a(\lambda)=\frac{1}{2\pi}\frac{\sqrt{(a_+-\lambda)(\lambda-a_-)}}{a\lambda}\1_{[a_-,a_+]}(x)\rmd x,\;a_\pm=(1\pm\sqrt{a})^2\;.
\end{equation}
In $F^{R_N}(x)$,  the $\lambda_k$'s denote the eigenvalues of $R_N$.
\end{lemma}

Our first main result is the $\sqrt{n}$-consistency of the estimator $\hat{\eta}$.

\begin{theorem}\label{th:consistency}
Let $\Y=(Y_1,\dots,Y_n)'$ satisfy Model (\ref{eq:modele_final}) with $\eta^{\star}>0$ and the entries $W_{i,j}$ of $\W$ satisfy Assumption \ref{assum:A1}. 
Then, for  all $q$ in  $(0,1]$, as $n,N\to\infty$ such that $n/N\to a \in (0,1]$,
$$
\sqrt{n}(\hat{\eta}-\eta^\star)=O_P(1).
$$
\end{theorem}

Such a result is a theoretical cornerstone to legitimate the use of an estimator. However, statistical inference has to be based on confidence sets.
The next step is thus to find the asymptotic distribution of $\sqrt{n}(\hat{\eta}-\eta^\star)$. 
Define for any $\eta\in[0,1]$ and $\lambda\geq 0$
$$
g(\eta,\lambda)=\frac{\lambda-1}{\eta (\lambda-1)+1}\;.
$$
Define also
$$
\gamma_n^{2}=\left\{\frac{1}{n}\sum_{i=1}^{n}g(\hat{\eta},\lambda_i)^{2}-\left(\frac{1}{n}\sum_{i=1}^{n}g(\hat{\eta},\lambda_i)\right)^{2}\right\}
$$
and
\begin{equation}\label{eq:sigma_tilde}
\gamma^2(a,\eta^\star)=\left\{\int g(\eta,\lambda)^2\rmd\mu_a(\lambda)
-\left(\int g(\eta,\lambda)\rmd\mu_a(\lambda)\right)^2\right\}\;.
\end{equation}


We are now ready to state our second main result about the asymptotic distribution of  $\sqrt{n}(\hat{\eta}-\eta^\star)$.
For general $q$, the result only holds when the entries of $\Z$, that is the random variables $Z_{i,j}$ are i.i.d. standard Gaussian. Indeed, as may be seen when computing the variances, we need to be able to find the asymptotic behavior of empirical means of functions of the eigenvalues together with the eigenvectors of the matrix $\R=\Z\Z'/N$.

\begin{theorem}\label{th:modele_general}
Let $\Y=(Y_1,\dots,Y_n)'$ satisfy Model (\ref{eq:modele_final}) with $\eta^{\star}>0$ and assume that the random variables $Z_{i,j}$ are i.i.d. ${\mathcal N}(0,1)$.
Then for any $q\in(0,1]$, as $n,N\to\infty$ such that $n/N\to a >0$,
$$\sqrt{n}(\hat{\eta}-\eta^\star)$$
 converges in distribution to a centered Gaussian random variable with variance
$$
\tau^2(a,\eta^\star, q)=\frac{2 }{\gamma^2(a,\eta^\star)} + 3 \frac{a^{2}{\eta^{\star}}^{2}}{\gamma^4(a,\eta^\star)}\left(\frac{1}{q}-1\right)S(a,\eta^\star)$$
where
$$
S(a,\eta^\star)=
\left[\int \frac{\lambda(\lambda-1)}{(\eta^\star(\lambda-1)+1)^{2}}\rmd\mu_a(\lambda)-\int\frac{\lambda}{(\eta^\star(\lambda-1)+1)}\rmd\mu_a(\lambda)\int\frac{\lambda-1}{(\eta^\star(\lambda-1)+1)}\rmd\mu_a(\lambda)\right]^{2}.
$$
\end{theorem}

In the case where $q=1$, the result holds in the general situation described in Assumption \ref{assum:A1}, and allows us
to propose confidence sets with precise asymptotic confidence level.
\begin{theorem}\label{th:modelenonsparse}
Let $\Y=(Y_1,\dots,Y_n)'$ satisfy Model (\ref{eq:modele_final}) with $q=1$ and with  $\eta^{\star}>0$. Assume also that the entries $W_{i,j}$ of $\W$ satisfy Assumption \ref{assum:A1}
then, as $n,N\to\infty$ such that $n/N\to a >0$,
$$\gamma_n\sqrt{\frac{n}{2}}\left(\hat{\eta}-\eta^\star\right)
$$ 
converges in distribution to ${\mathcal N}(0,1)$.
\end{theorem}

Let us now give some additional comments on the previous results.
Firstly, it has to be noticed that none of the limiting variance depends on $\sigma^\star$. Secondly,
Theorem \ref{th:modele_general} is proved here only in the case where the $Z_{i,j}$ are i.i.d. Gaussian. This is because we used several times that the matrix of eigenvectors of $\Z\Z'/N$ is independent of the eigenvalues, and uniformly distributed on the set of orthonormal matrices. 
We think that the result of  Theorem \ref{th:modele_general} is also valid when the $Z_{i,j}$ are defined from the $W_{i,j}$ satisfying Assumption \ref{assum:A1}, 
as suggested by the numerical results obtained in Section \ref{sec:numeric}.
To prove it requires  new results in an active research topic of the random matrix theory field.
We can observe in the expression of $\tau^2(a,\eta^\star)$ given in Theorem \ref{th:modele_general}
that the presence of $q$ is counterbalanced by the presence of $a^2$. This will be confirmed by the results
obtained in the numerical results given in Section \ref{sec:numeric}. Finally, we can observe that $2/(n\gamma_n^2)$ 
corresponds to the usual inverse of the Fisher information associated to $\eta$. This result is classical in the case where 
$N$ is fixed and $n$ tends to infinity but did not exist in the framework where both $n$ and $N$ tend to infinity even if
it was already used in biological applied papers for deriving standard errors and confidence intervals. 
Theorem \ref{th:modelenonsparse} proves that this result still holds even in the case
where both $n$ and $N$ tend to infinity.

To the best of our knowledge, the effect of the presence of null components in the random effects has never been taken into account for computing the asymptotic variance of an estimator of the heritability. 
This is the contribution of Theorem \ref{th:modele_general}. This theorem shows that the asymptotic variance contains 
an additional term which increases its value in the case $q<1$ with respect to the case $q=1$. It is shown in Section 3.3
how the computation of the asymptotic variance can be altered if this additional term is neglected.
In practical situations, computing the standard error given by Theorem \ref{th:modele_general}
requires the knowledge of $q$ which is in general unknown. However, if an estimation of $q$ is available for any practical
reasons, the result of Theorem  \ref{th:modele_general} can be used for computing confidence intervals and standard errors, 
see Section \ref{sec:discussion} for further details.

%


\section{Numerical experiments}\label{sec:numeric}

In this section, we first explain how to implement our method and then we illustrate the theoretical results 
of Section \ref{sec:theory} on finite sample size observations for both cases: $q=1$ and $q<1$. 
We also compare the results obtained with our approach
to those obtained by the GCTA software described in \cite{Yang:2010} and
\cite{yang:lee:goddard:visscher:2011} which is a reference 
in quantitative genetics.

\subsection{Implementation}

\begin{figure}[!ht]
\begin{center}
\includegraphics[width=0.6\textwidth]{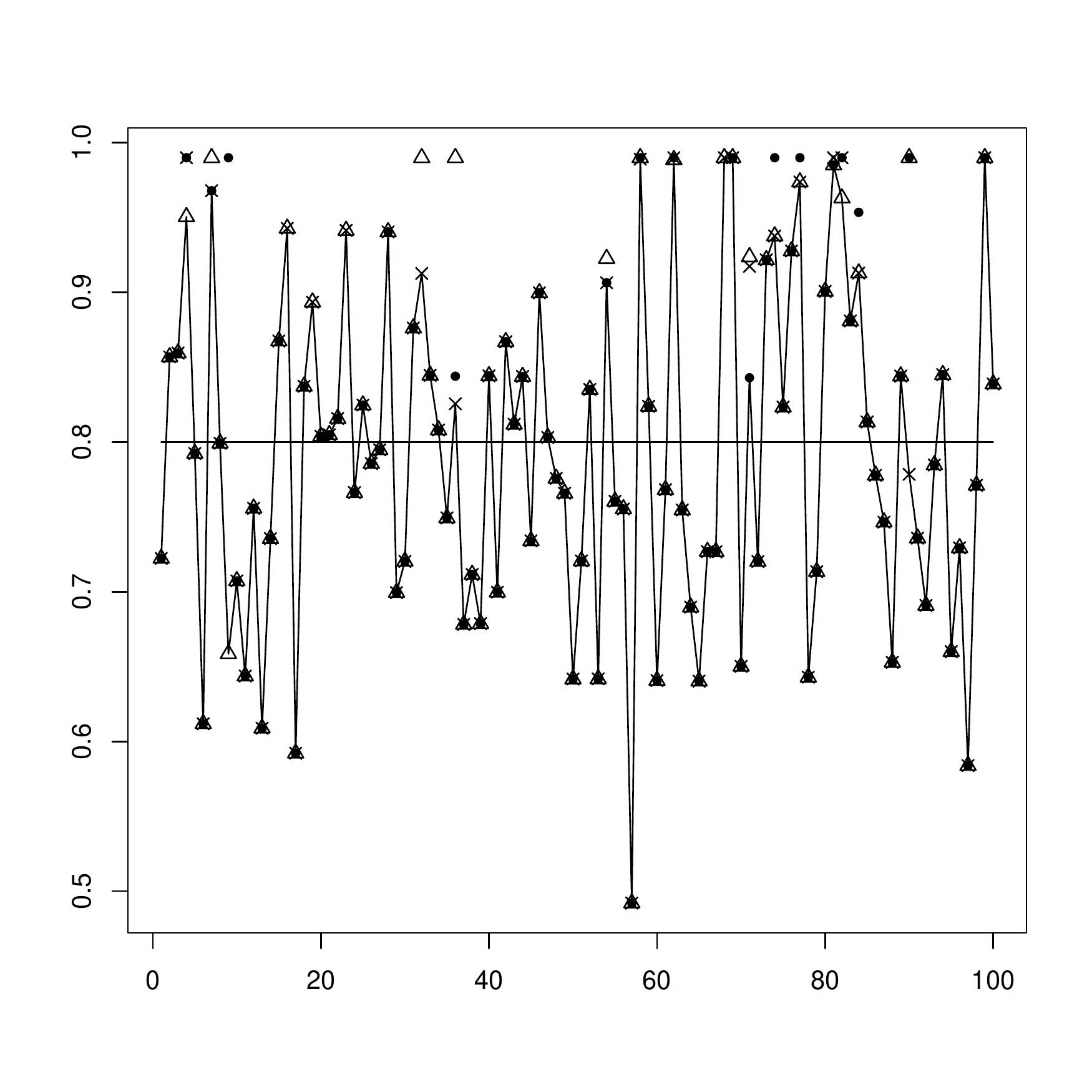} 
\end{center}
\caption{Estimation of $\hat{\eta}$ obtained in the case $a=0.1$ and $\eta^\star=0.8$ for different values of initialization: $\eta^{(0)}=0.1$ (dots), $\eta^{(0)}=0.5$ (triangles) and $\eta^{(0)}=0.9$ (crosses). The plain line displays the estimations obtained with 
our method to select the best initialization value and the $x$-axis is the replication number.}
\label{fig:init}
\end{figure}

In order to obtain $\hat{\eta}$, we used a Newton-Raphson approach which is based on the following
recursion: starting from an initial value $\eta^{(0)}$,
$$
\eta^{(k+1)}= \eta^{(k)} - \frac{L_n'(\eta^{(k)})}{L_n''(\eta^{(k)})}\;,\; k\geq 1\;,
$$ 
where $L_n'$ and $L_n''$ denote the first and second derivatives of $L_n$ defined in (\ref{eq:Ln}), respectively.
The closed form expression of these quantities are given in (\ref{eq:Ln_prime}) and (\ref{eq:Ln_second}), respectively. 
In practice, this approach converges after at most 20 iterations and is not very sensitive to the initialization, namely to
the value of $\eta^{(0)}$.  However, in particular cases, the value of the initialization can have an influence 
on the estimation of $\eta^\star$. This is the case, for instance, when the real value $\eta^\star$ is close to 1. 
In these situations, our algorithm can provide an estimation bigger than 1 and we constrained our method to return a 
value equal to 0.99. Figure \ref{fig:init} shows the estimations obtained on $100$ replications when $a=0.1$ and $\eta^\star=0.8$. 
From this figure, we can see that the estimation of $\eta^\star$ does not depend in general on the initialization, 
except in some cases. Moreover, the best choice for $\eta^{(0)}$ is not constant from one replication to another.  
In order to limit the effect of the initialization, our algorithm uses several values for  $\eta^{(0)}$ and 
whenever the estimations differ, it keeps the estimation which is the farthest away from the boundaries.

\subsection{Results in Model (\ref{eq:distrib_u_e_comp_nulles}) when $q=1$}

\begin{figure}[!ht]
\begin{center}
\begin{tabular}{ccc}
\hspace{-5mm}\includegraphics[width=40mm]{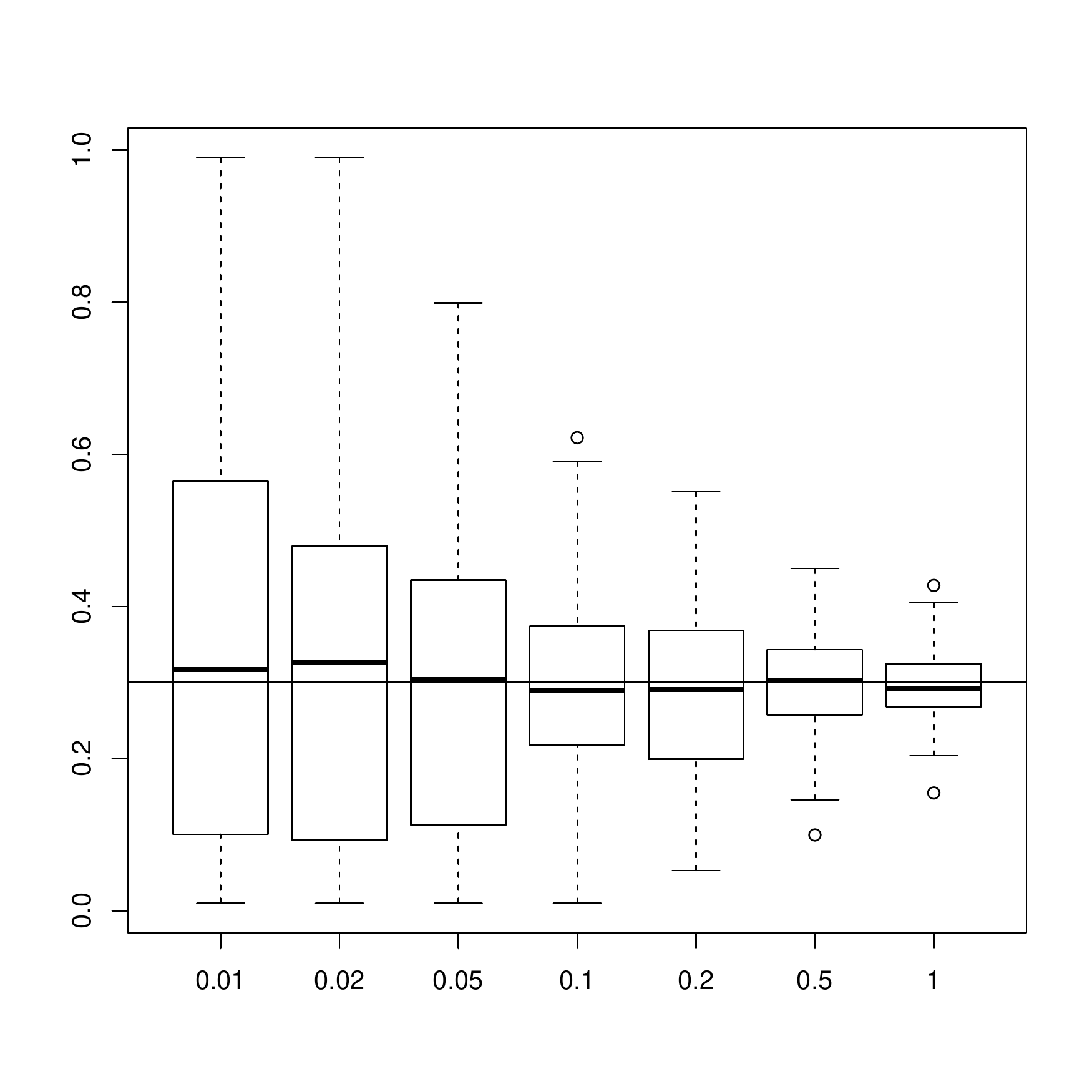} &
\includegraphics[width=40mm]{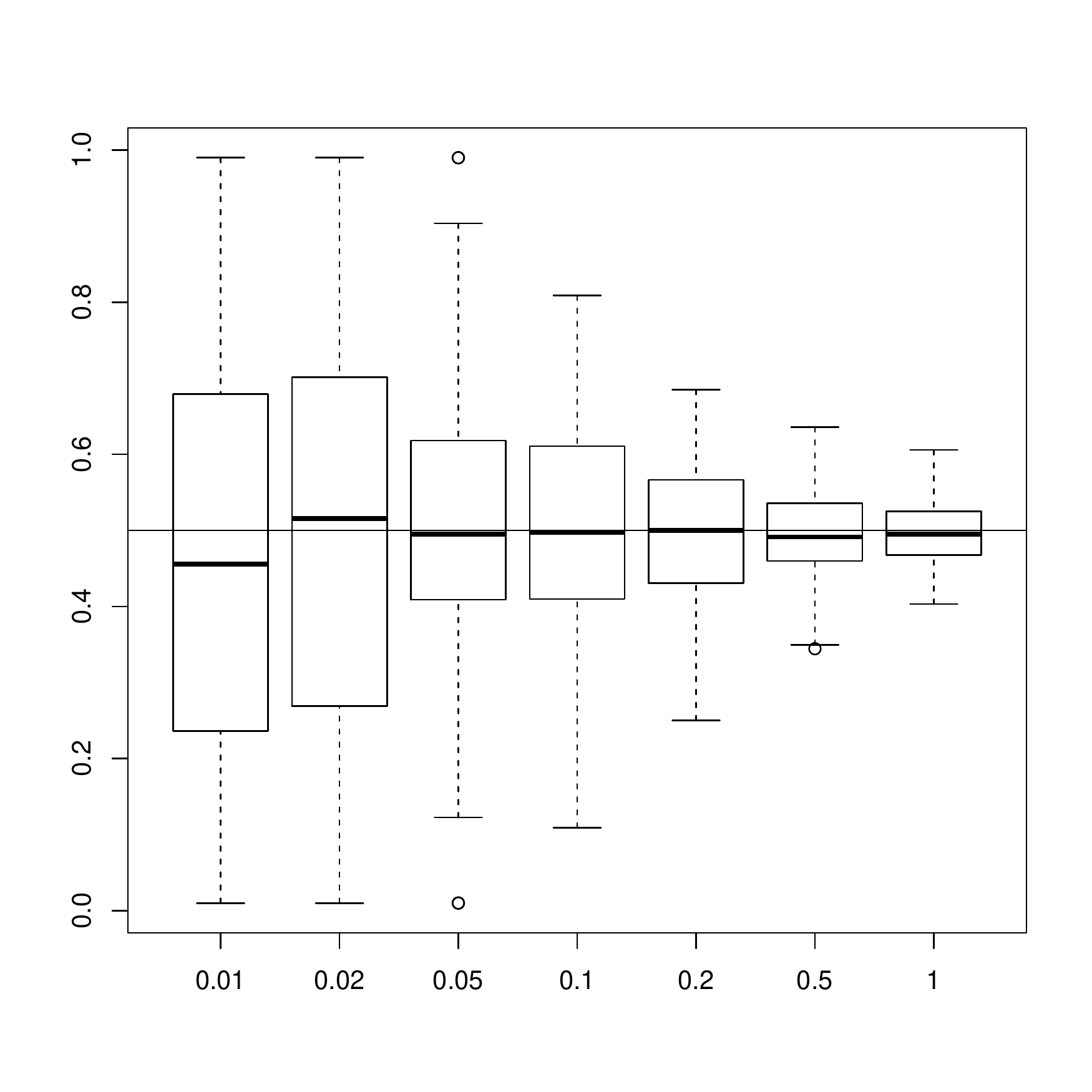} &
\includegraphics[width=40mm]{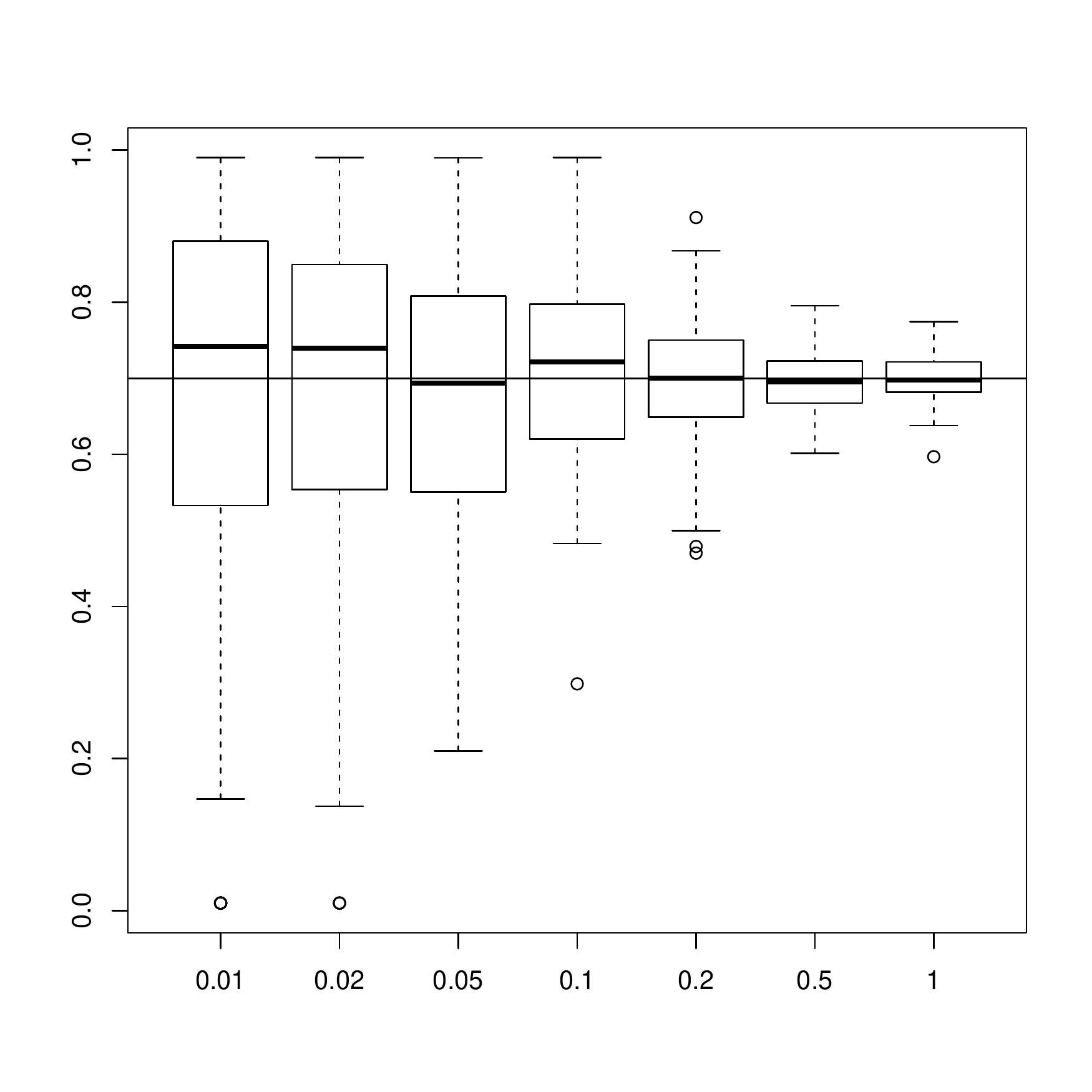}
\end{tabular}
\end{center}
\caption{Boxplots of $\hat{\eta}$ for different values of $a$, for $\eta^{\star}=0.3$ (left), $\eta^{\star}=0.5$ (middle) and $\eta^{\star}=0.7$ (right). The horizontal line corresponds to the true value of $\eta^\star$. The whiskers of each boxplot correspond to the first and third quartiles.}
\label{fig:boxplot_eta}
\end{figure}

 We shall first consider the performance of the estimator $\hat{\eta}$ when $q=1$
for $\eta^\star$ in $\{0.3, 0.5, 0.7\}$, $n=1000$, $\sigma_u^\star=0.1$ and 
for $a$ in  $\{0.01,0.02,0.05,0.1,0.2,0.5,1\}$, where $a=n/N$. We generated 500 data sets according to
Model (\ref{eq:modele}) using these parameters and $\Z$ as defined in (\ref{eq:normalization_1})
where the $W_{i,j}$ are binomial random variables with parameters 2 and $p_j$. In our experiments the $p_j$'s 
are uniformly drawn in $[0.1,0.5]$.
The corresponding boxplots of $\hat{\eta}$ are displayed
in Figure \ref{fig:boxplot_eta}. We can see from this figure that our approach provides
unbiased estimators of $\eta^\star$ and that the smaller the $a$ the larger the empirical variance. 

\begin{figure}[!ht]
\begin{center}
\begin{tabular}{ccc}
\hspace{-5mm}\includegraphics[width=40mm]{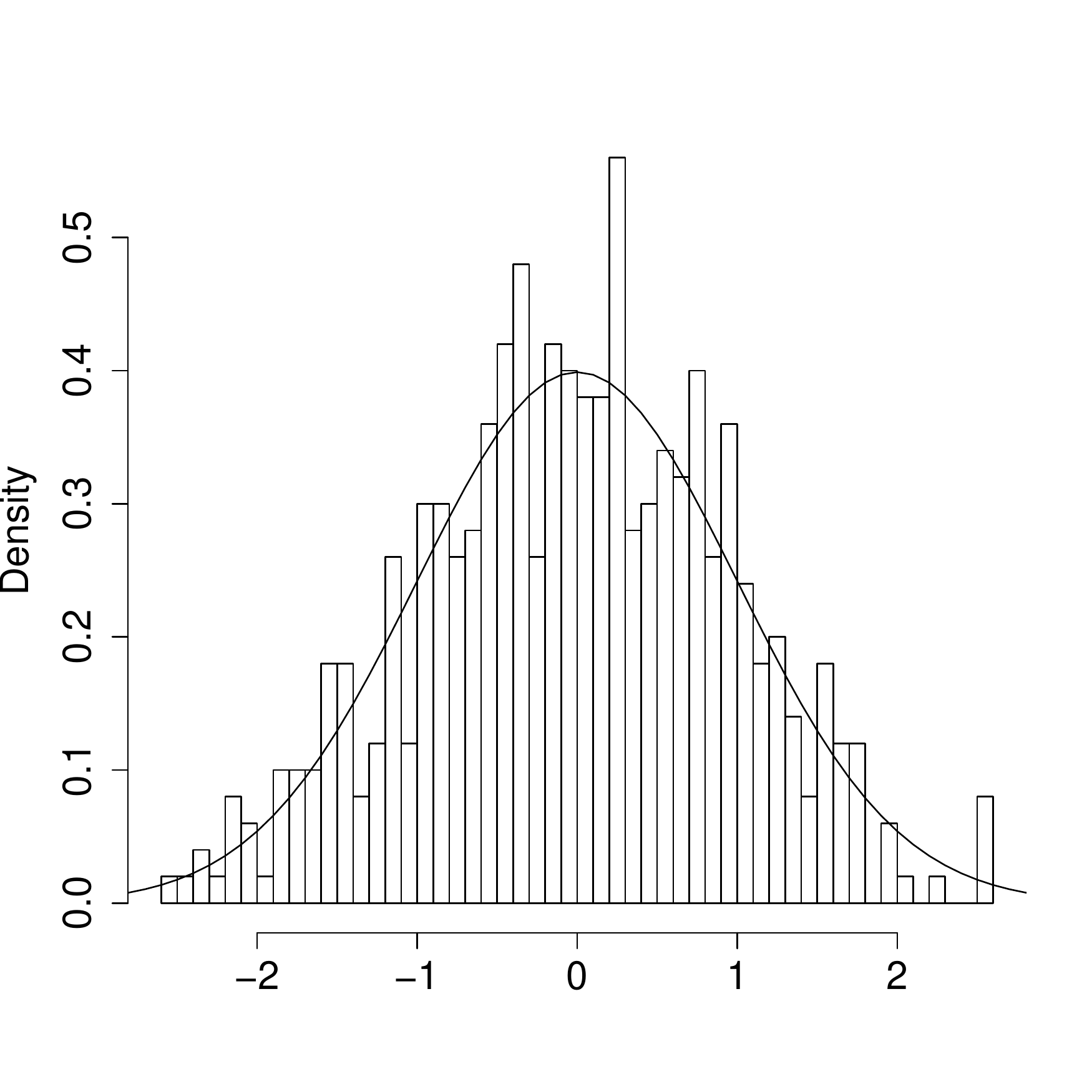} &
\includegraphics[width=40mm]{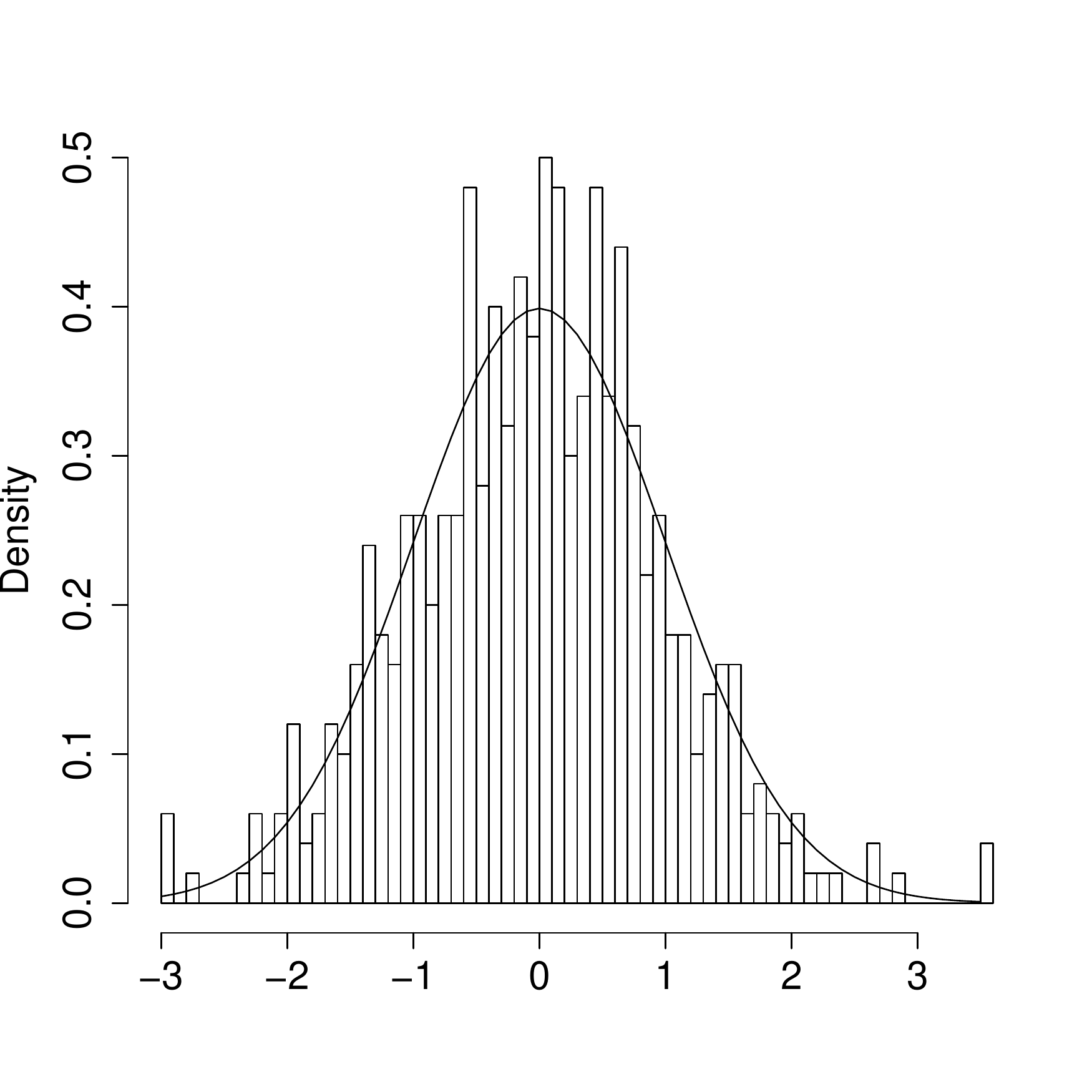} &
\includegraphics[width=40mm]{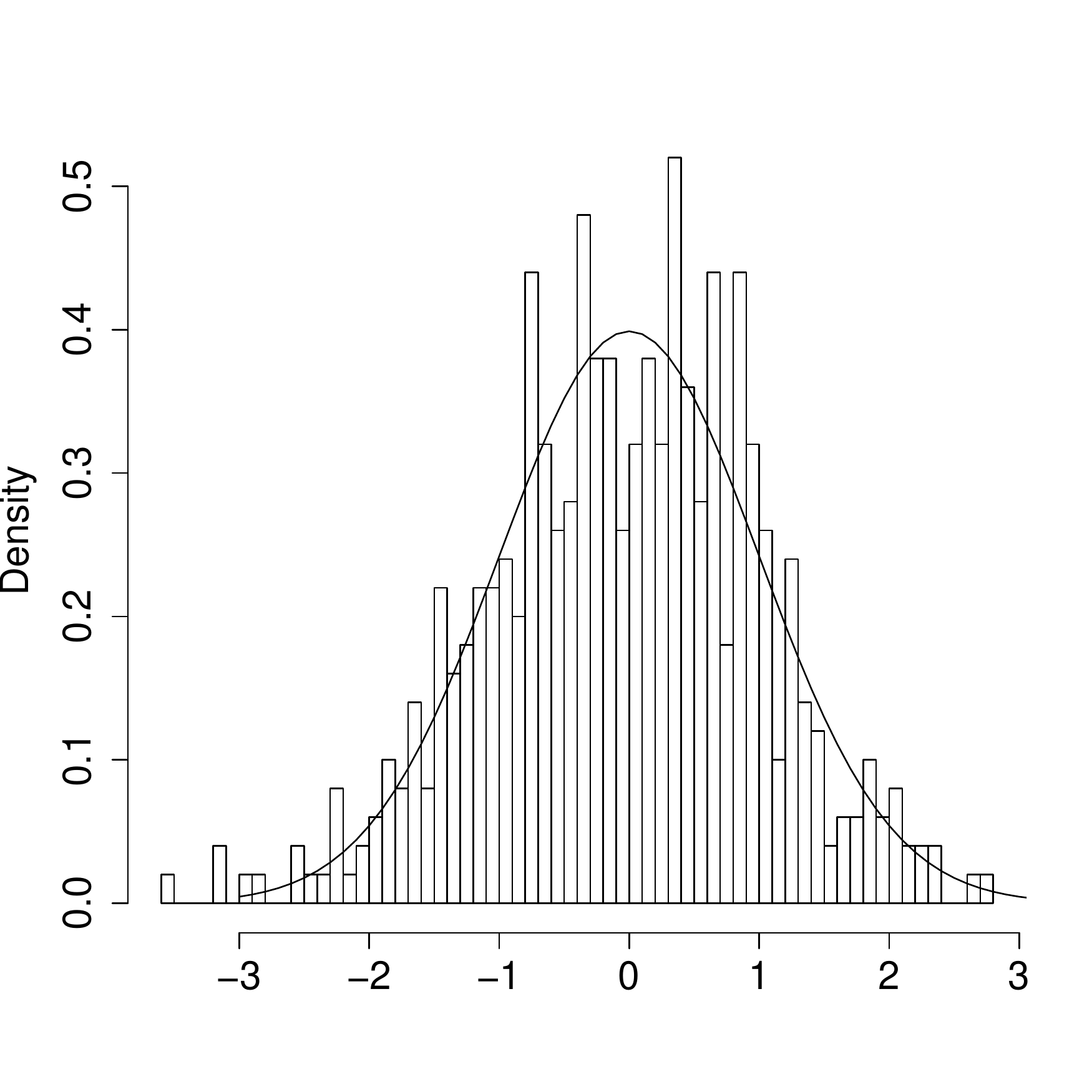}
\end{tabular}
\caption{Histograms of $\gamma_{n}(n/2)^{1/2}\left(\hat{\eta}-\eta^\star\right)$ for $\eta^\star=0.5$ and $a=0.05$ (left), 
$a=0.1$ (middle), $a=0.5$ (right) and the p.d.f of a standard Gaussian random variable in plain line.}
\label{fig:hist_eta07_q1}
   \end{center}
\end{figure}

\begin{figure}[!ht]
\begin{center}
\begin{tabular}{cc}
\includegraphics[width=0.35\textwidth]{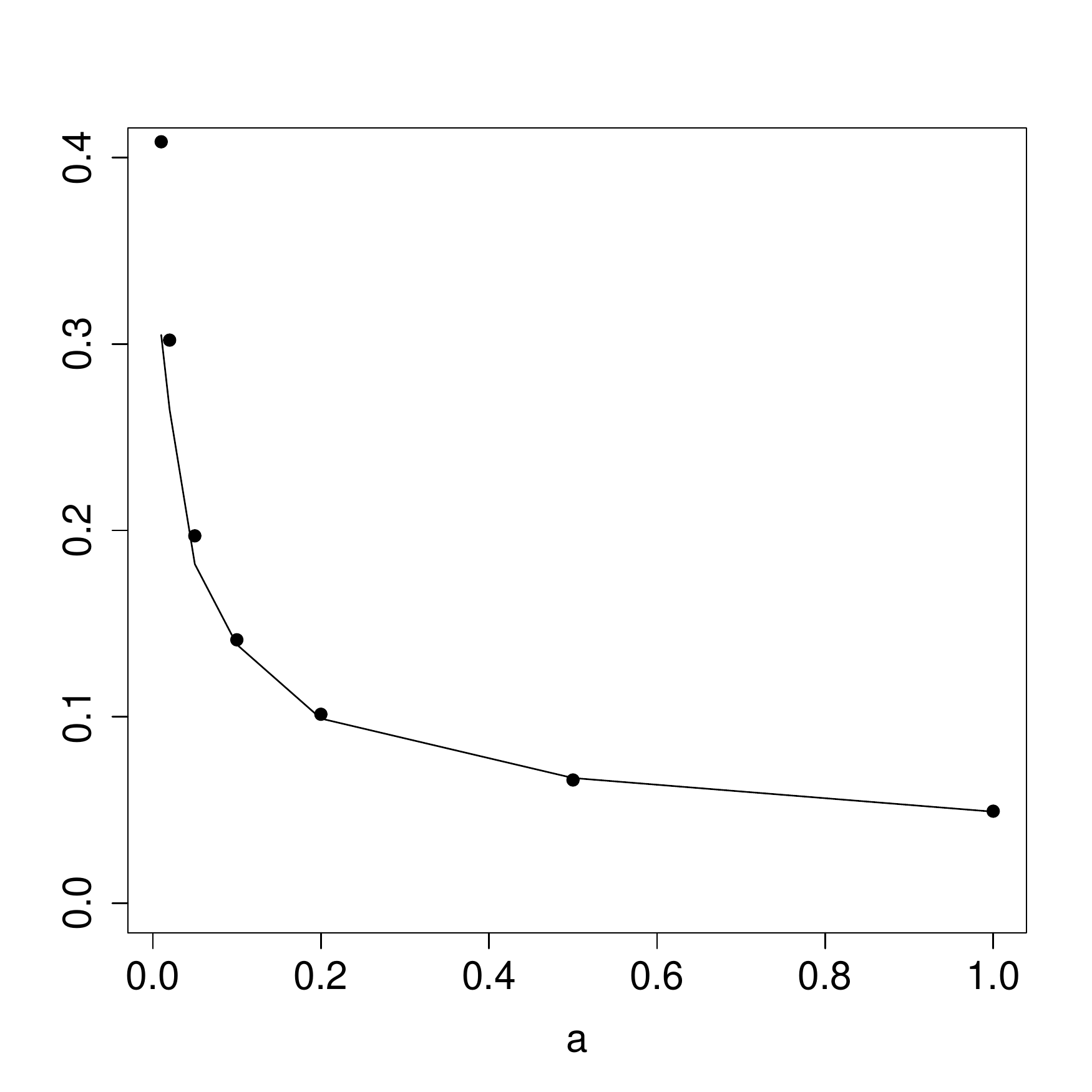} &
\includegraphics[width=0.35\textwidth]{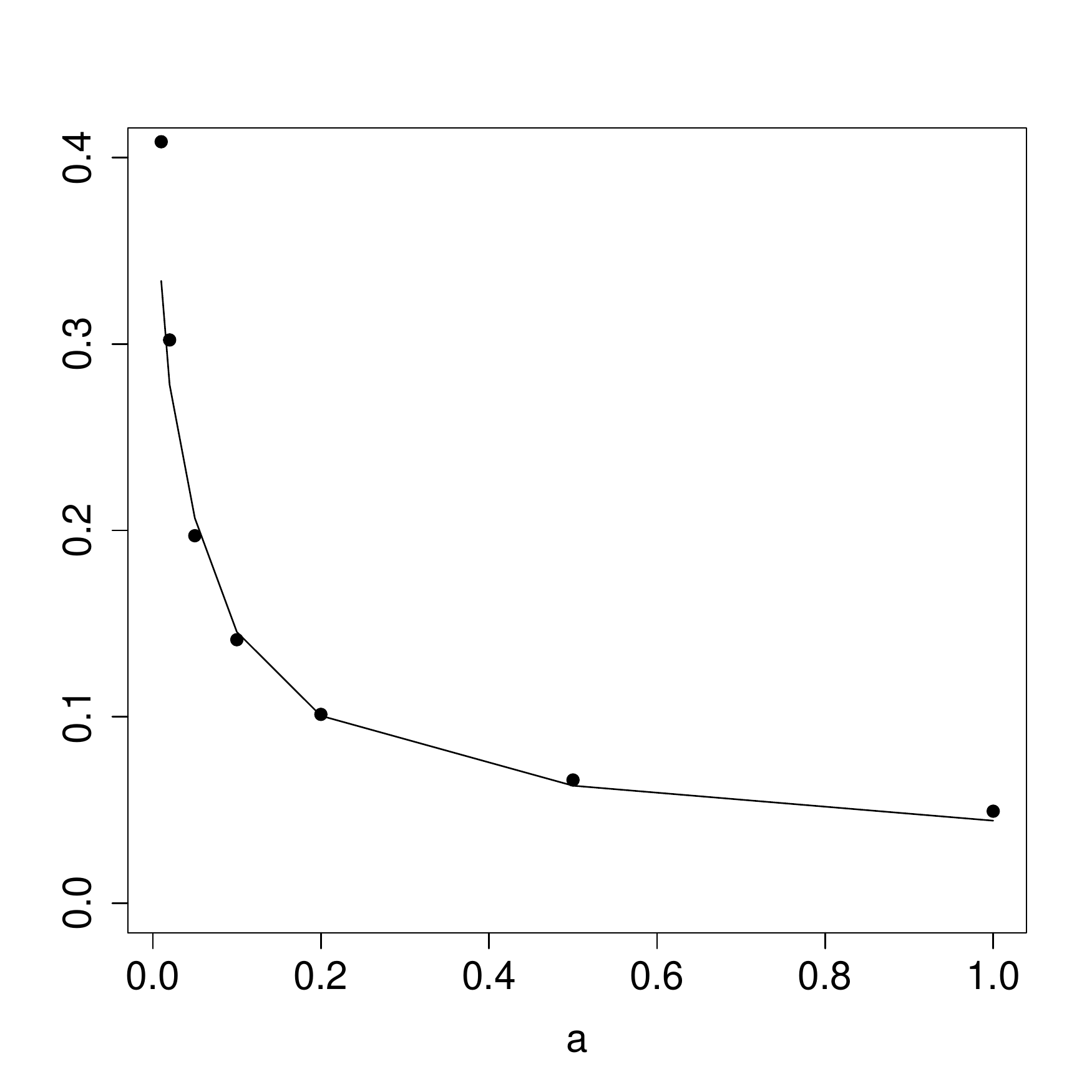} 
\end{tabular}
\end{center}
\caption{Values of $n^{-1/2}\sqrt{2\gamma_n^{-2}}$ (``$\bullet$'')  and the empirical standard deviation of 
$(\hat{\eta}-\eta^{\star})$ (plain line) for several values of $\eta^{\star}$ (0.3 (left), 0.5 (right)).
}
\label{fig:compar_var}
\end{figure}

In order to illustrate the central limit theorem given in Theorem \ref{th:modelenonsparse}, we first display
in Figure \ref{fig:hist_eta07_q1} the histograms of $\gamma_{n}(n/2)^{1/2}\left(\hat{\eta}-\eta^\star\right)$
along with the p.d.f of a standard Gaussian random variable for $\eta^{\star}=0.5$ and different values of $a$.
We can see that the Gaussian p.d.f fits well the data in all the considered cases.
We also display in Figure \ref{fig:compar_var} the values of $n^{-1/2}\sqrt{2\gamma_n^{-2}}$ and the empirical standard deviation
of $(\hat{\eta}-\eta^{\star})$ averaged over all the experiments. As shown in Theorem \ref{th:modelenonsparse}, we also observe empirically that
both quantities are very close.

In practice, the value of $\gamma_{n}^{-1}(n/2)^{-1/2}$ can be used for deriving confidence intervals
for $\eta^{\star}$. As we can see from Figure \ref{fig:compar_var}, our approach leads to
very accurate confidence intervals for $a$ larger than 0.1 even in finite sample size cases.

\begin{figure}[!ht]
\begin{center}
\includegraphics[width=0.5\textwidth]{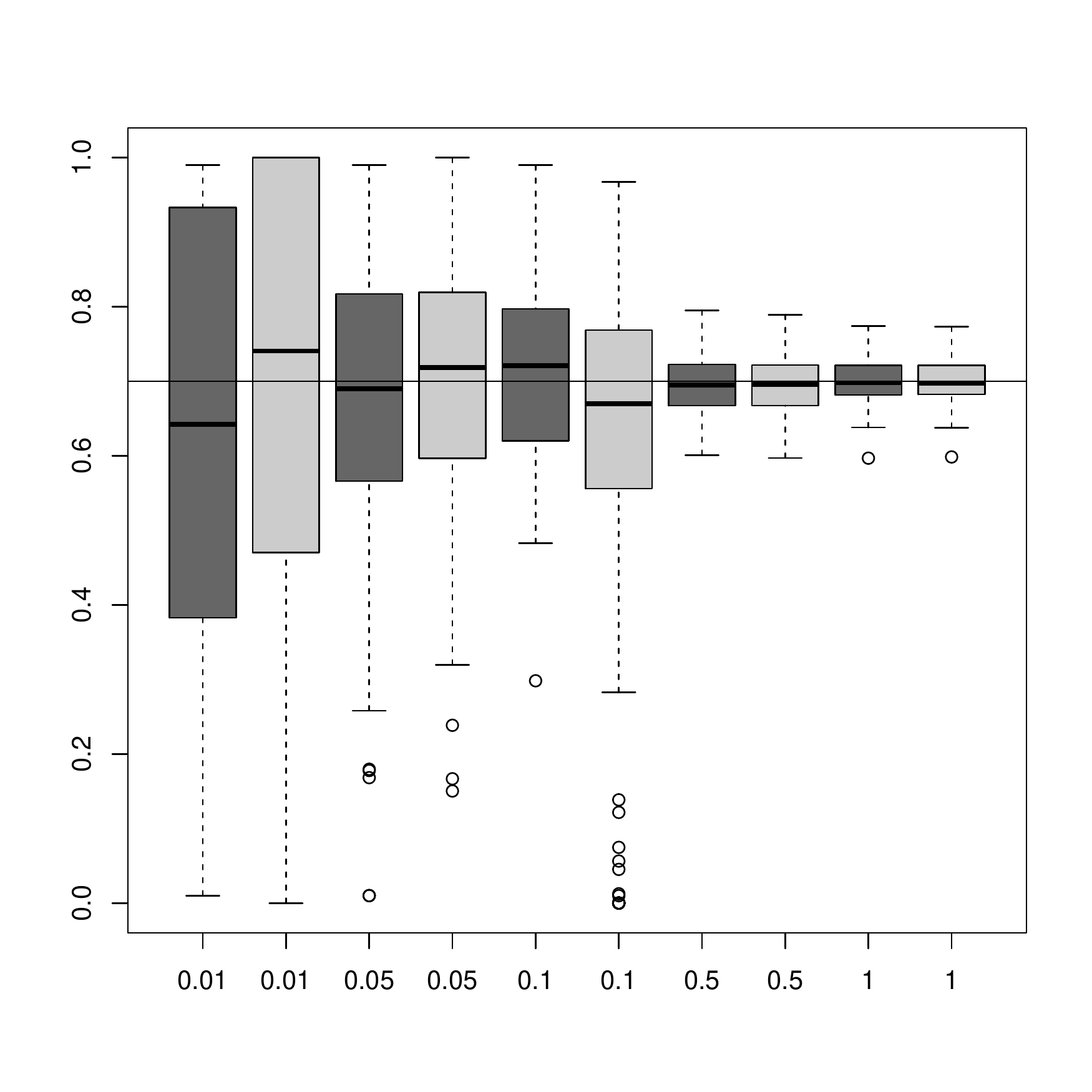} 
\end{center}
\caption{Boxplots of $\hat{\eta}$ for different values of $a$, using our method (dark gray) and GCTA (light gray). The whiskers of each boxplot are the first and third quartiles.}
\label{fig:compar_GCTA}
\end{figure}

Let us now compare our results with those obtained with the software GCTA.
As we can see from Figure \ref{fig:compar_GCTA} which displays the boxplots of $\hat{\eta}$ for different values of $a$
when $\eta^\star=0.7$, the results found by our approach and GCTA are very close. In both cases, we observe that
when $a$ is close to $1$ the estimations of $\eta^\star$ are very accurate but when $a$ is small the standard error 
becomes very high.

 \subsection{Results in model \ref{eq:distrib_u_e_comp_nulles} when $q<1$}





This section is dedicated to the study of the performance of $\hat{\eta}$ when $q<1$.
We generated 500 data sets according to
Model (\ref{eq:modele}) for $\eta^{\star}=0.7$, $a \in \{0.05,0.1,0.5,1 \}$, different values of $q$ and $\Z$ defined in (\ref{eq:normalization_1})
where the $W_{i,j}$ are binomial random variables with parameters 2 and $p_j$. In our experiments the $p_j$'s 
are uniformly drawn in $[0.1,0.5]$.

\begin{figure}[!ht]
\begin{center}
\begin{tabular}{cccc}
\includegraphics[width=50mm]{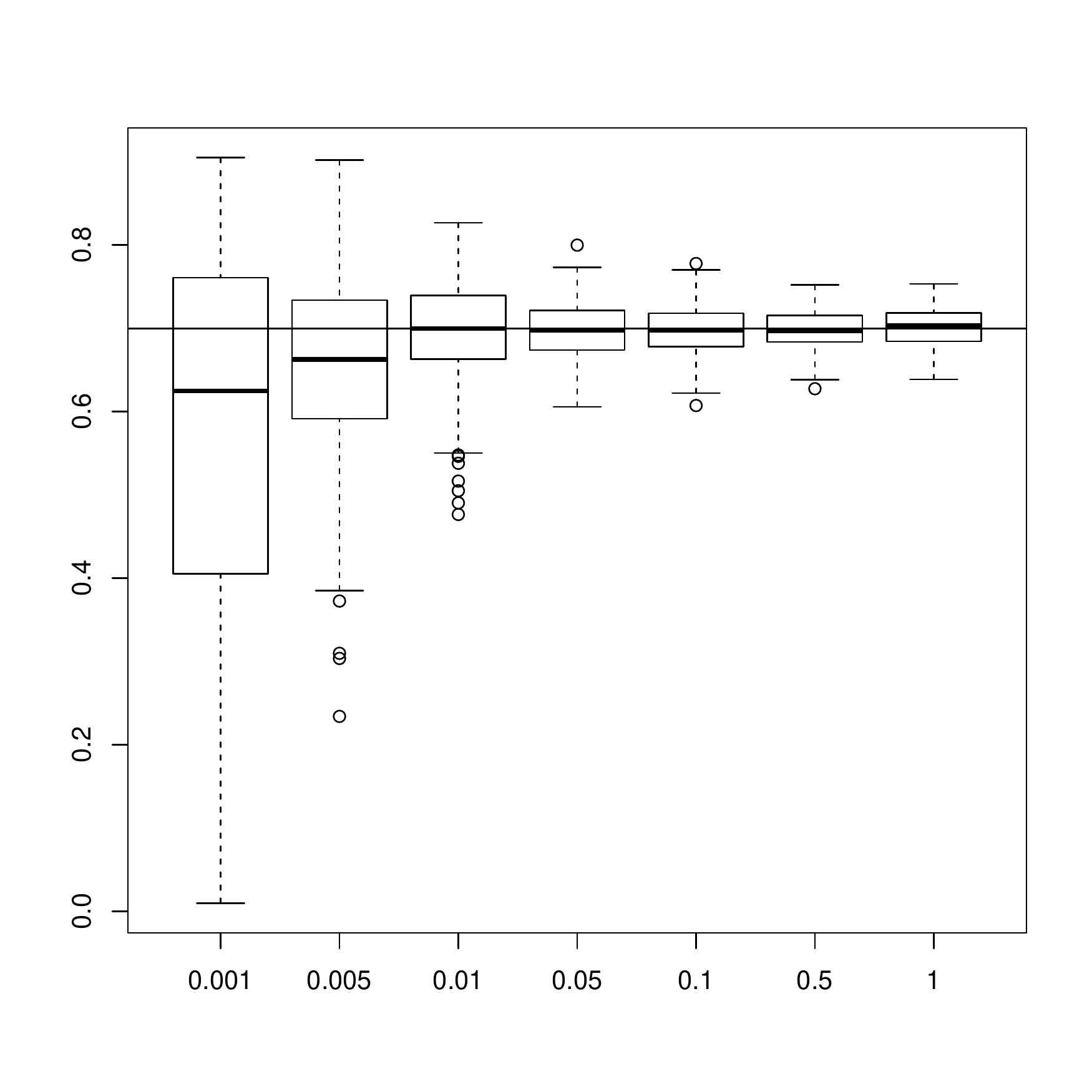} &
\includegraphics[width=50mm]{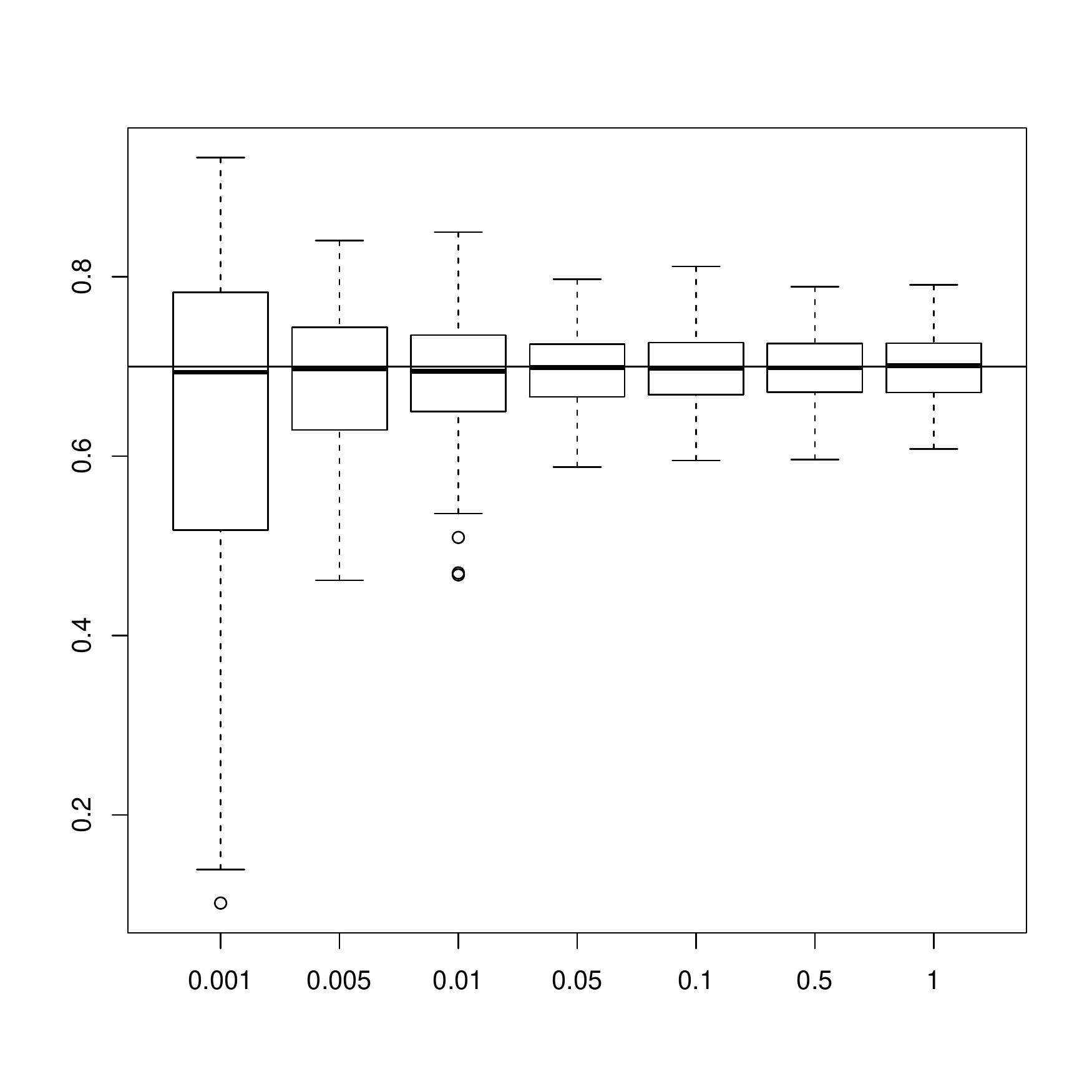} &\\
\includegraphics[width=50mm]{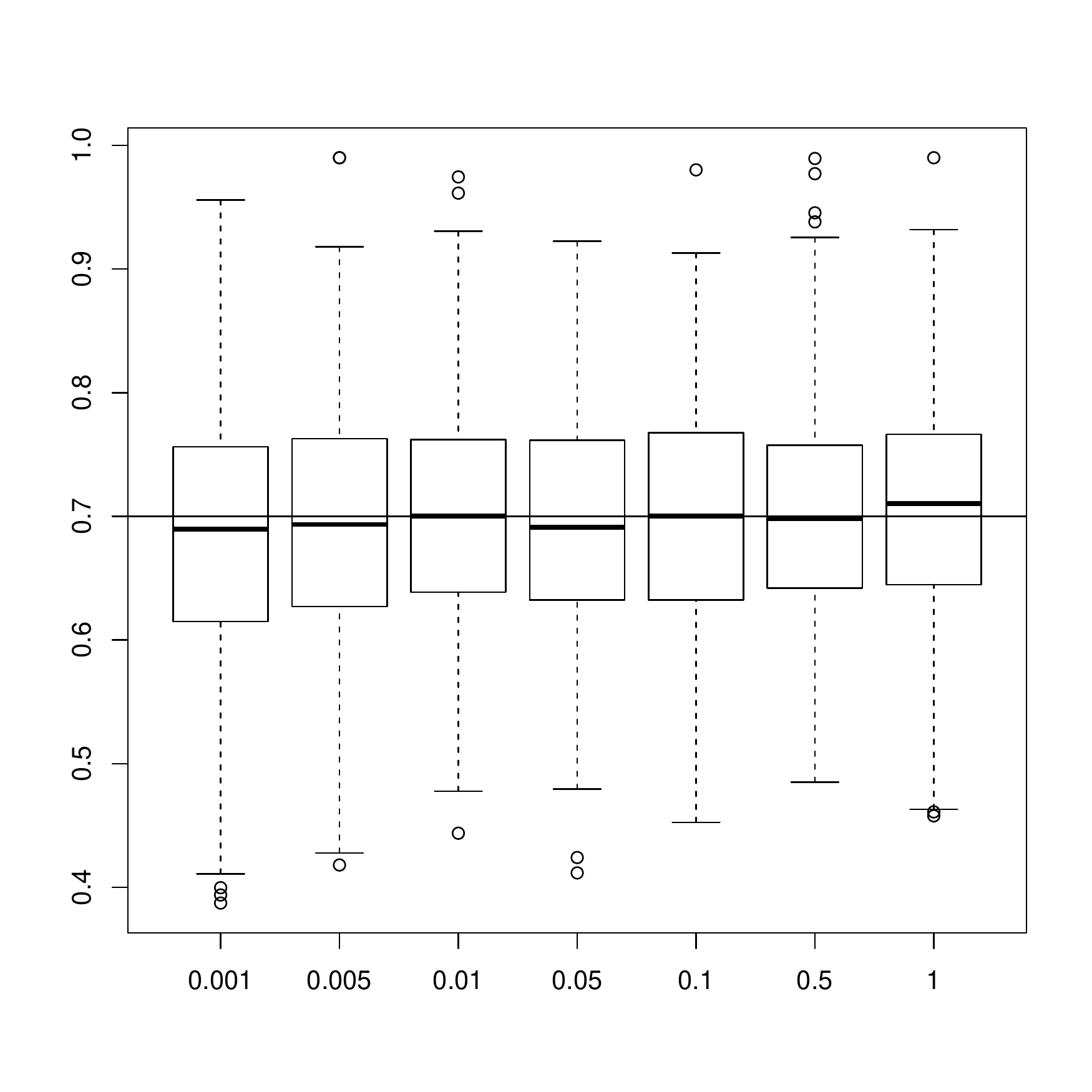} &
\includegraphics[width=50mm]{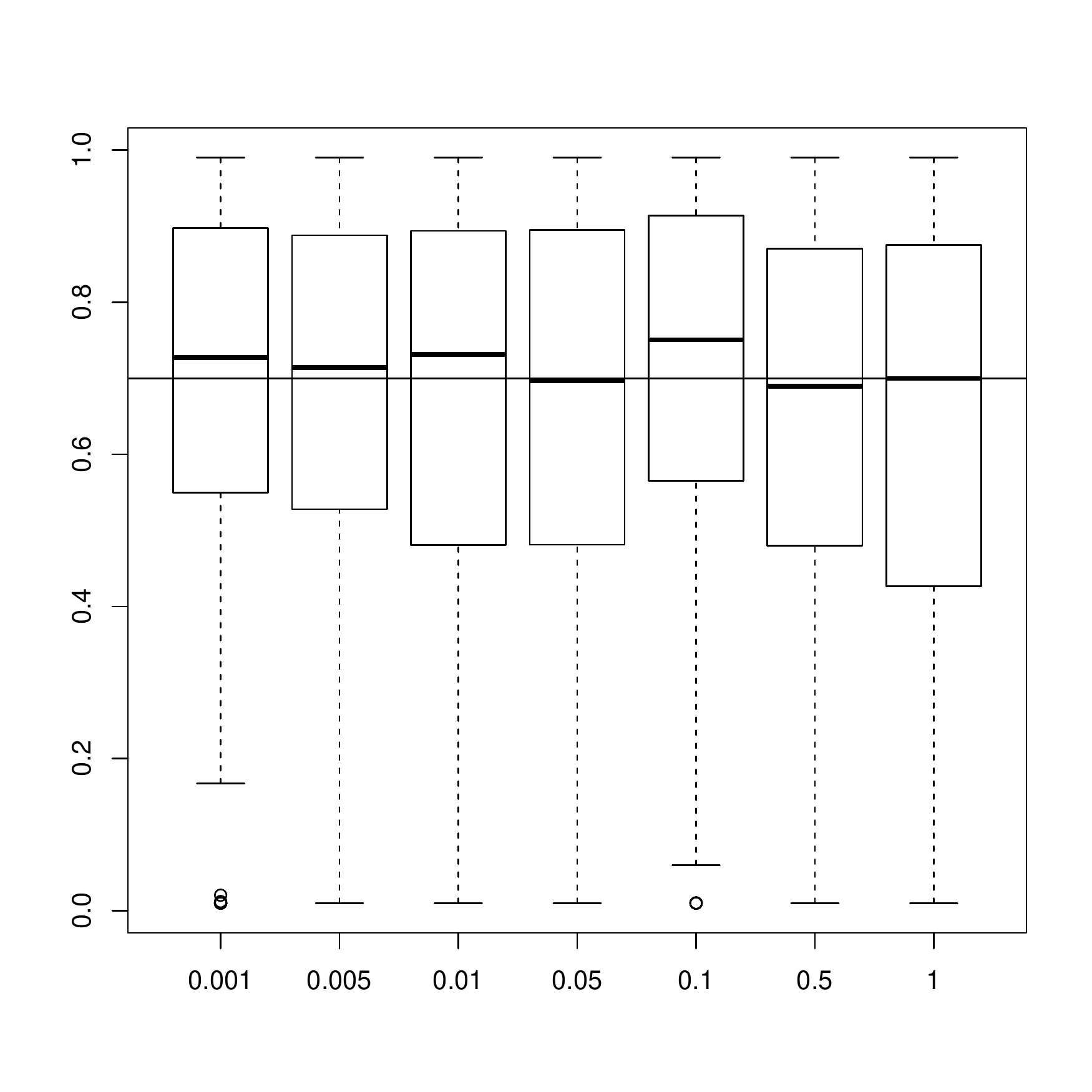}
\end{tabular}
\caption{Boxplots of $\hat{\eta}$ for different values of $q$, with $\eta^{\star}=0.7$ and $a=1$ (top left), $a=0.5$ (top right), $a=0.1$ (bottom left) and $a=0.01$ (bottom right). 
The boxplots are based on $500$ replications. The whiskers of each boxplot are the fist and third quartile.}
\label{fig:boxplot_etahat_compnulles}
   \end{center}
\end{figure}

Figure \ref{fig:boxplot_etahat_compnulles} displays the boxplots of $\hat{\eta}$ for these parameters. 
We can see from this figure that for small values of $a$, the estimators of $\eta^{\star}$
have the same behavior for $q=1$ and $q<1$. However, when $a=1$ or $a=0.5$, we can see from this figure
that the presence of null components strongly alter the performance of the estimator of $\eta^\star$.
Since in typical GWAS experiments, $a=0.01$ or even smaller, the results of Figure \ref{fig:boxplot_etahat_compnulles}
could lead to conclude that considering the case $q<1$ is not necessary for such values of the parameter $a$. 
However, as already noticed from Figure \ref{fig:boxplot_eta}, the variance of $\hat{\eta}$ is very large
for small values of $a$, hence considering the presence of null components and proposing a strategy for selecting
only the non null components of $\u$ could be one way to increase $a$  and thus to diminish the variance of
$\hat{\eta}$.

\begin{figure}[!ht]
\begin{center}
\begin{tabular}{cc}
\includegraphics[width=0.4\textwidth]{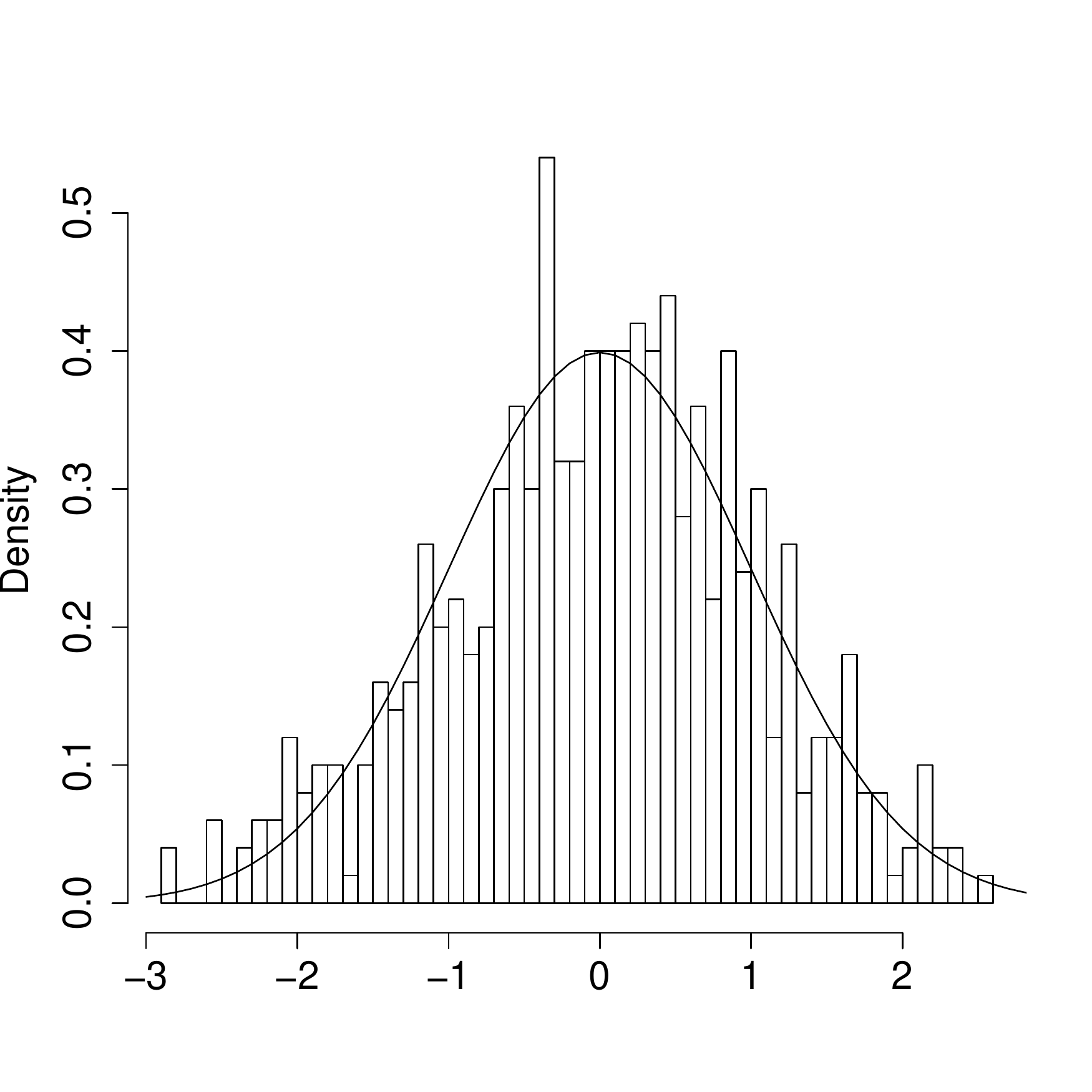}&\includegraphics[width=0.4\textwidth]{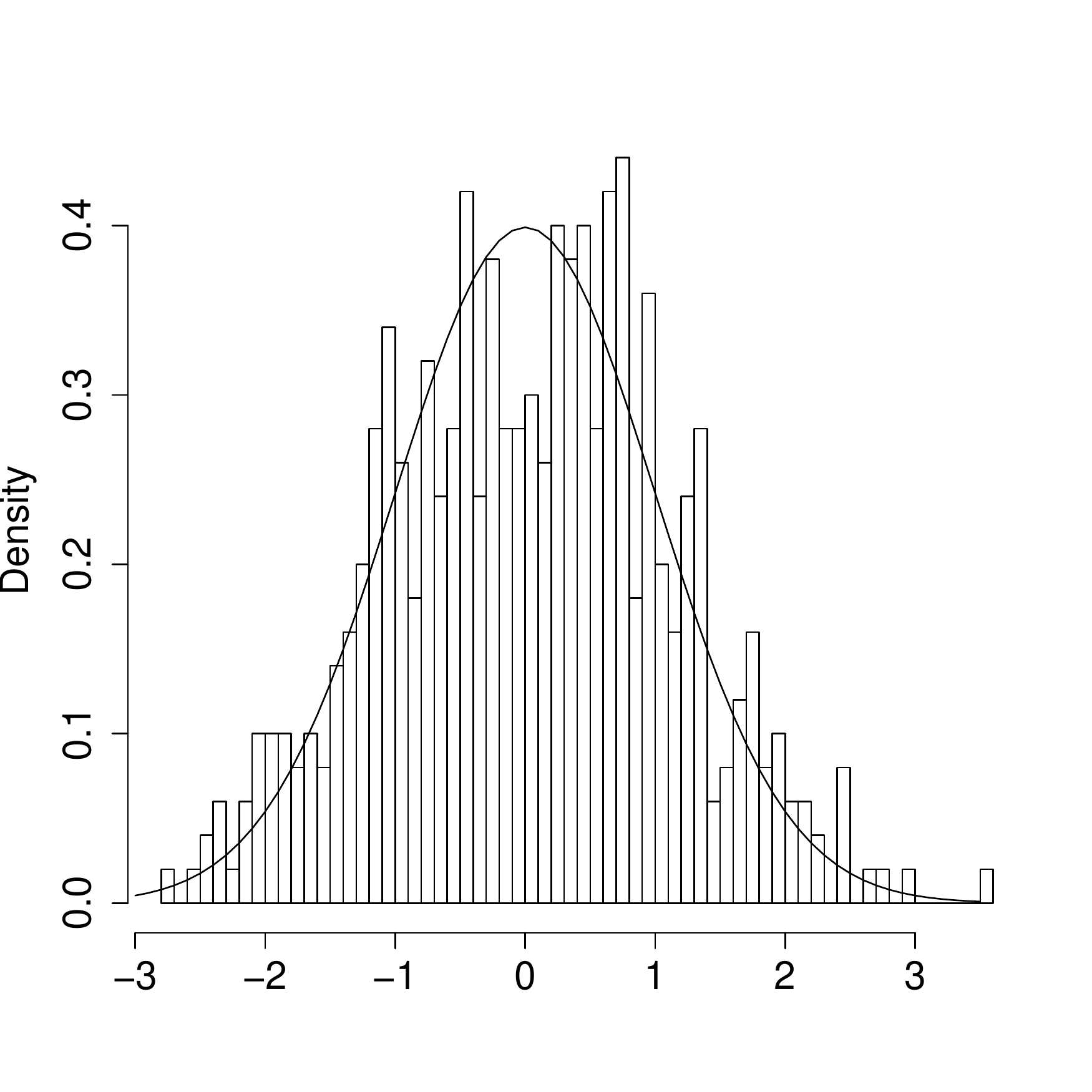}\\
\includegraphics[width=0.4\textwidth]{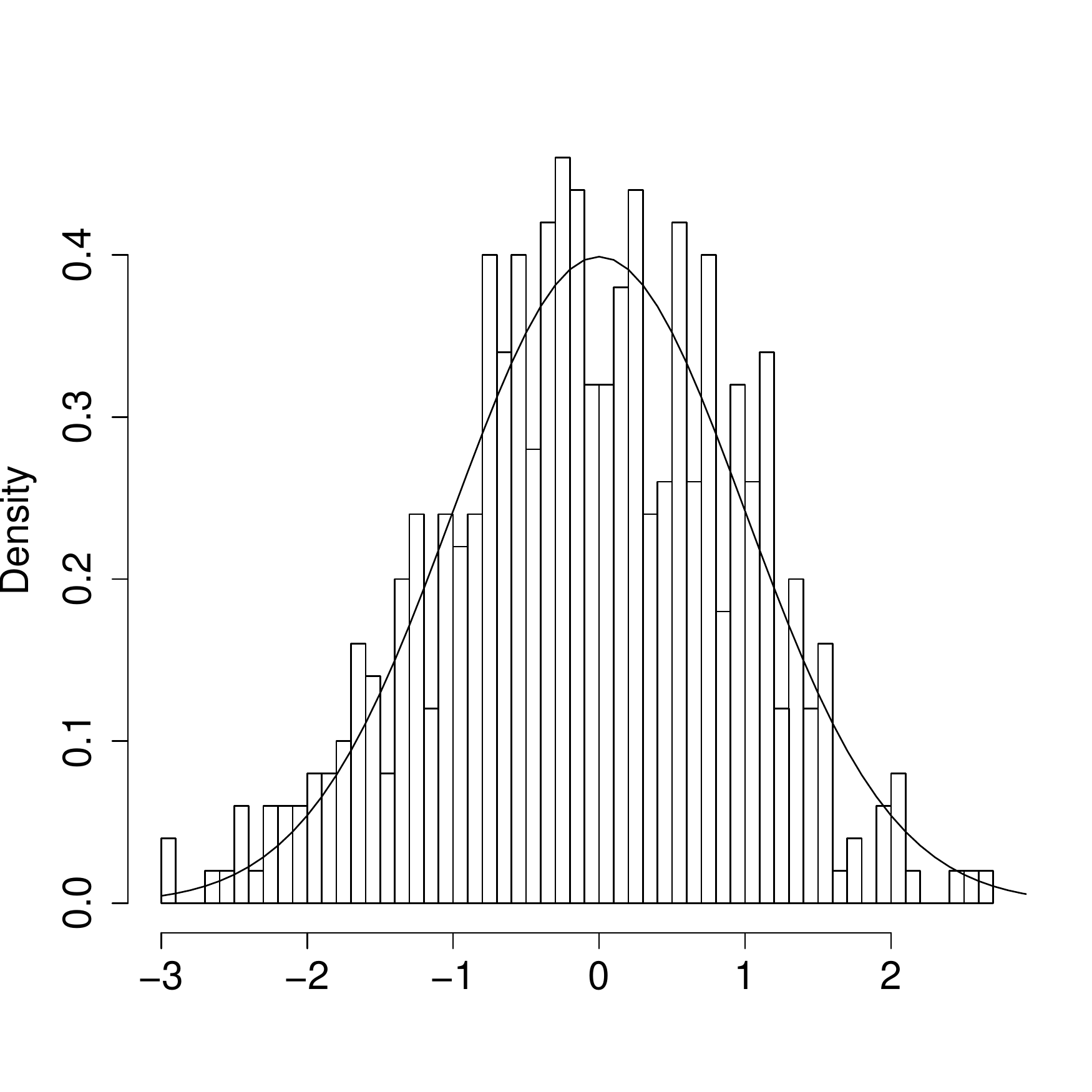}&\includegraphics[width=0.4\textwidth]{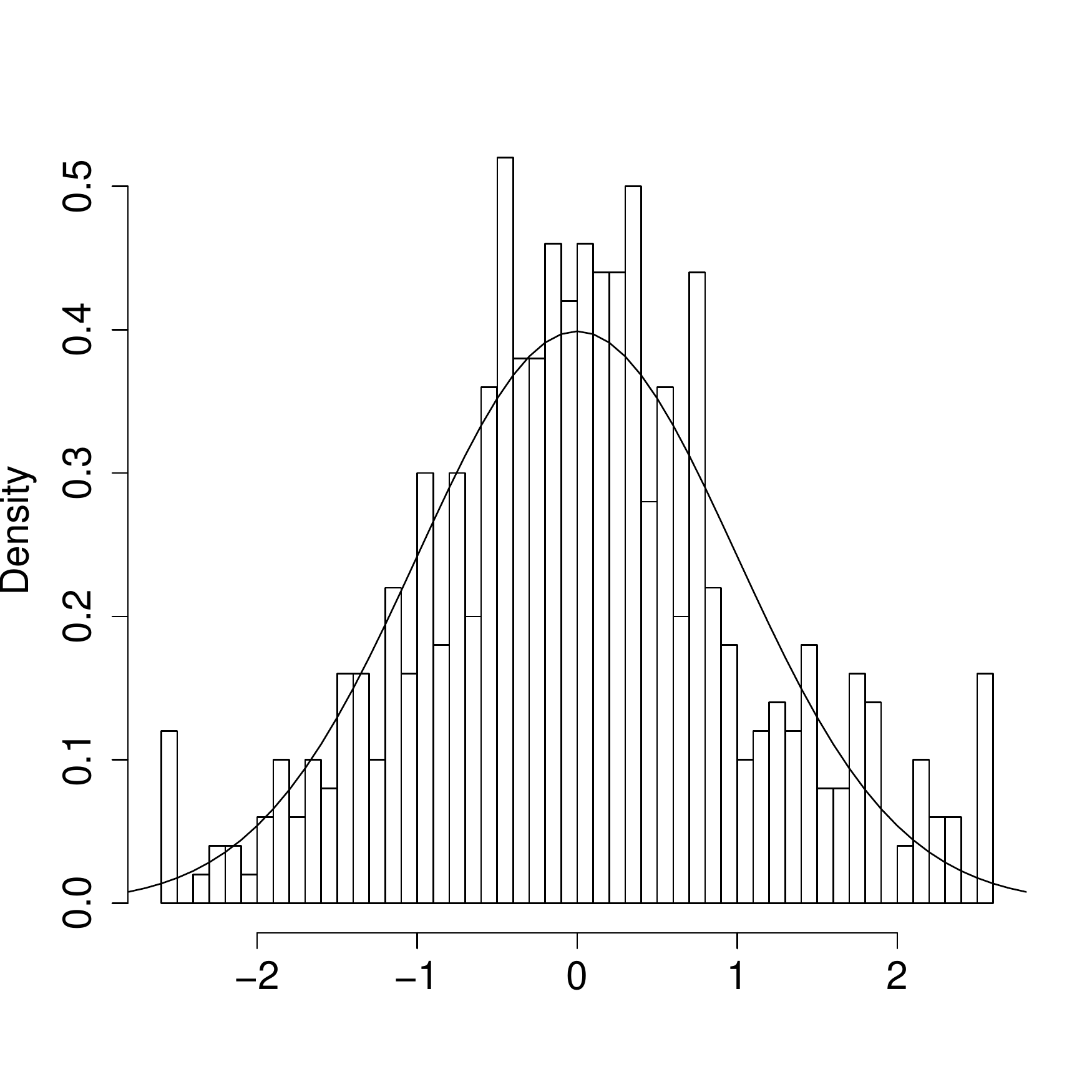}
\end{tabular}
\caption{Histograms of $\tau_{n}^{-1} n^{1/2}\left(\hat{\eta}-\eta^\star\right)$ for $a=0.5$ and $q=0.5$ (top left), $a=0.1$ and $q=0.1$ (top right),
and for $a=0.1$ and $q=0.01$ (bottom left), $a=0.05$ and $q=0.1$ (bottom right).}
\label{fig:hist_eta07_differ_a_differ_q}
\end{center}
\end{figure}

In order to illustrate the central limit theorem given in Theorem \ref{th:modele_general}, we first display
in Figure \ref{fig:hist_eta07_differ_a_differ_q} the histograms of $\tau_{n}^{-1} n^{1/2}\left(\hat{\eta}-\eta^\star\right)$
along with the p.d.f of a standard Gaussian random variable for $\eta^{\star}=0.7$, two values of $q$: $q=0.01$ and $q=0.1$ and $a=0.5$ (top)
and two values of $a$: $a=0.2$ and $a=0.5$ with $q=0.5$ (bottom).
Here, $\tau_{n}$ is the empirical version of $\tau(a,\eta^\star,q)$  where $\gamma$ is replaced by $\gamma_n$ and
$S(a,\eta^\star)$ is replaced by its empirical version with the eigenvalues of $\R$.
When $a$ is large ($a=0.5$), we can see that the higher $q$ the better the Gaussian p.d.f fits the histograms.


\begin{figure}[!ht]
\begin{center}
\includegraphics[angle=-90,width=0.6\textwidth]{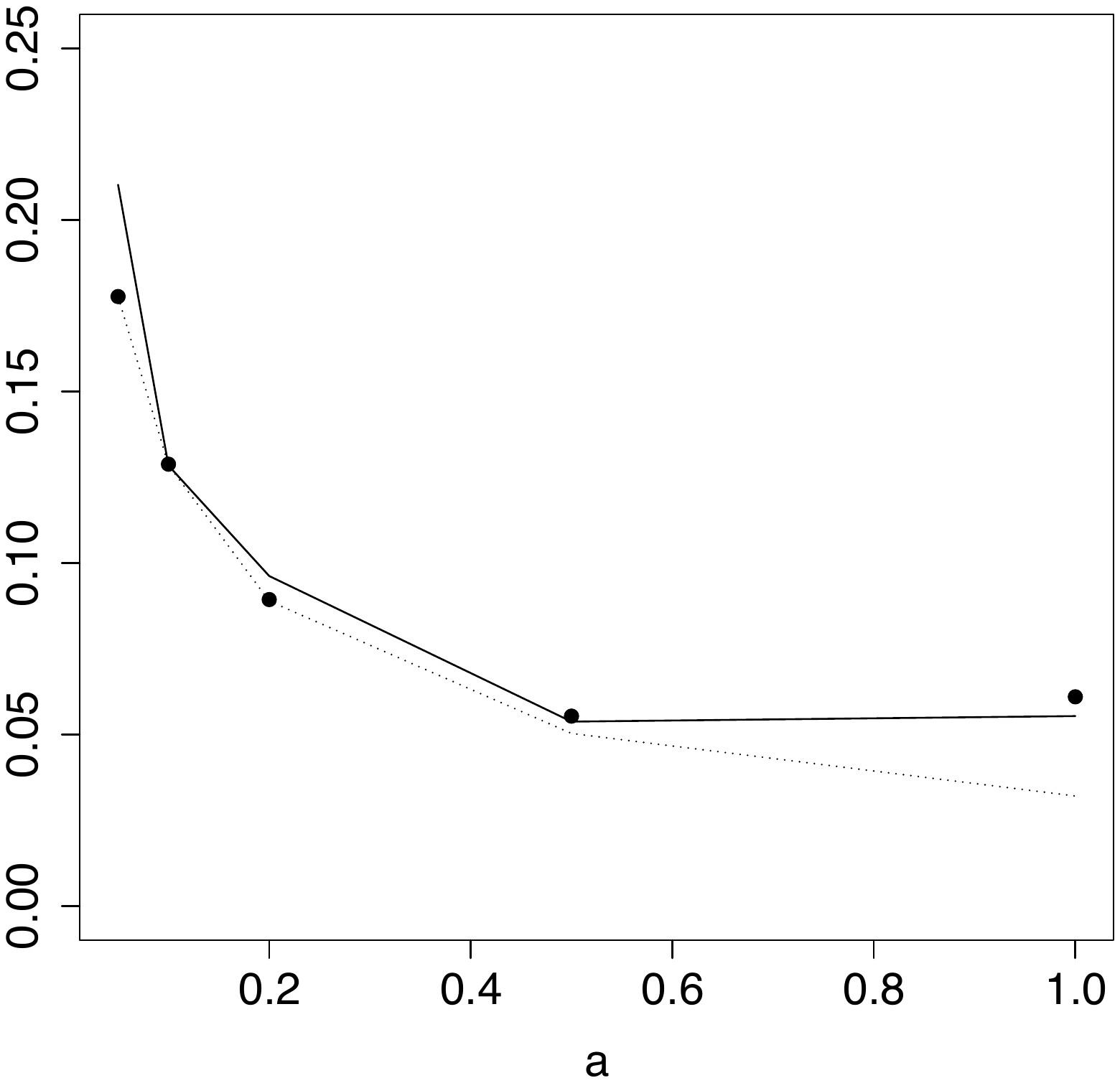} 
\end{center}
\caption{Values of $n^{-1/2}\tau_{n}$ with the real value of $q$ ($q=0.5$) (``$\bullet$''), $q=1$ (dotted line)
and the empirical standard deviation of $(\hat{\eta}-\eta^{\star})$ (plain line) for $\eta^{\star}=0.7$.}
\label{fig:compar_var_q}
\end{figure}

We also display in Figure \ref{fig:compar_var_q} the values of $n^{-1/2}\tau_{n}$ and the empirical standard deviation
of $(\hat{\eta}-\eta^{\star})$ averaged over all the experiments for $\eta^\star=0.7$ and $q=0.5$. 
As shown in Theorem \ref{th:modele_general}, we observe empirically that both quantities are very close.
We also display in this figure the value of $n^{-1/2}\tau_{n}$ with $q=1$ which boils down to consider the asymptotic standard deviation
found in the non sparse model. We can see from this figure that neglecting the term depending on $q$ leads to underestimate
the asymptotic variance of $\hat{\eta}$ and that this difference is all the more striking that $a$ is close to 1.

\section{Discussion}\label{sec:discussion}

In the course of this study, we have proposed a methodology for estimating the 
heritability in high dimensional linear mixed models. This methodology has two main features. Firstly, 
the theoretical performances of our estimator are established in a non standard theoretical framework where
$n$ and $N$ tend to infinity and where the components of the random effect part can be equal to zero.
Secondly, the computational burden of our approach is very low which makes its use possible on real data coming from GWAS
experiments.

As a byproduct of the central limit theorem that we establish for $\eta^\star$ we can derive 
a confidence interval for the heritability. However, the confidence intervals depend on $q$ which is
the proportion of non null components in $\u$ and which is general unknown. 
For estimating $q$, several strategies can be considered. One could, for instance, use a GWAS approach 
to compute the $p$-values of the correlation tests of each SNP with the observations $\Y$ and then keep only the most significant ones.
Such a practical approach can be used for providing a lower bound for $q$. A refinement of this approach has been proposed by \cite{toro:2014}
who observed, through numerical studies, that for a fixed value of the heritability, the minimal $p$-value is all the more low 
that the number of causal SNPs is small. Hence, performing a GWAS approach on a given data set allows them to have an idea of the number
of SNPs which are likely to be causal.
One could also propose another practical method based on a variable selection technique. 
Such an approach has already been proposed by \cite{Fan_Li} in the context
of sparse linear mixed models. However, the framework in which their theoretical results are derived is different from the one that
is considered in our paper. We are currently working on a paper \cite{papier_appli_nous:2014} which presents a variable selection
method which is adapted to our framework and which could be used for estimating the proportion $q$ of non null components in the random effects.
%
%

Moreover, we did not take into account the linkage disequilibrium issue which would require to extend our results to the case  where
the columns of the random matrix are correlated. This question will be the subject of a future work.

Let us write the singular value decomposition (SVD) of the $n\times N$ matrix $\Z/\sqrt{N}$ as
$$
\frac{1}{\sqrt{N}}\Z=\U\left(\sqrt{\D}\;\mathbf{0}\right)\V'
$$
where $\U$ (already introduced in Section \ref{sec:intro}) is a $n\times n$ orthonormal matrix, $\V$ is a $N\times N$ orthonormal matrix and
$\sqrt{\D}$ is a $n\times n$
diagonal matrix having its diagonal entries equal to $\sqrt{\lambda_i}$, the $\lambda_i$'s being the eigenvalues
of $\R=\Z\Z'/N$ previously defined. 
Thus, \eqref{eq:modele_final} rewrites as
\begin{equation}\label{eq:modele_svd}
\widetilde{\Y}=\U'\Y=\left(\sqrt{\D}\;\mathbf{0}\right)\V'\t+\sigma^\star\sqrt{1-\eta^\star}\,\widetilde{\eps}\;,
\end{equation}
where $\widetilde{\eps}=\U'\eps$ is a $n\times 1$ centered
Gaussian vector having a covariance matrix equal to identity.

We shall use repeatedly the following lemma which is proved in Section \ref{tech:lem}.

\begin{lemma}\label{lem:var_quad_form}
Let $\widetilde{\Y}$ be defined by (\ref{eq:modele_svd}) and $\H$ be a $n\times n$ diagonal matrix, then
$$
\Var\left(\widetilde{\Y}'\H\widetilde{\Y}|\Z\right)
=2{\sigma^\star}^4\Tr\left[\H^2\left\{(1-\eta^\star)\Id_{\rset^n}+\eta^\star\D\right\}^2\right]
+3{\sigma^\star}^4{\eta^\star}^2\left(\frac{1}{q}-1\right)\sum_{1\leq i\leq N} M_{ii}^2\;,
$$
where
$$
\M=\V\left(\begin{array}{cc}\D\H &0\\
0&0
\end{array}\right)\V'\;,
$$
and
$$
\Var\left(\widetilde{\Y}'\H\widetilde{\Y}|\Z\right)
\leq 2{\sigma^\star}^4\Tr\left[\H^2\left\{(1-\eta^\star)\Id_{\rset^n}+\eta^\star\D\right\}^2\right]
+3{\sigma^\star}^4{\eta^\star}^2\left(\frac{1}{q}-1\right)\Tr[ \D^{2} \H^{2}].
$$
\end{lemma}

Another useful lemma will be the following.
\begin{lemma}\label{lem:trickMP}
Under Assumption \ref{assum:A1}, let $h:\rset^{+}\to\rset^{+}$ be such that there exist $\alpha>0$ and $C$ such that for all $n$,
$$
\PE\left(\frac{1}{n}\sum_{i=1}^{n}h(\lambda_{i})^{1+\alpha}\right)\leq C.
$$
Then
$$
\frac{1}{n}\sum_{i=1}^{n}h(\lambda_{i})=\int h(\lambda)d\nu_{a}(\lambda)+o_{p}(1).
$$
\end{lemma}
The proof of this lemma follows from the application of Lemma \ref{lem:MP_noniid} to the bounded function $h\1_{h\leq M}$, and the Markov inequality applied to the empirical mean of $h\1_{h>M}$.

\begin{lemma}\label{lem:lambda2}
Under Assumption \ref{assum:A1} let $n,N\to\infty$ be such that $n/N\to a>0$. Then there exists $C$ such that for all $n$,
$$
\PE\left[\frac{1}{n}\sum_{i=1}^{n}\lambda_{i}^{2}\right]\leq C.
$$
\end{lemma}
To prove the lemma, notice that $\sum_{i=1}^{n}\lambda_{i}^{2}=\Tr[\Z\Z'/N^{2} ]$. But
\begin{eqnarray*}
\PE\left(\Tr\left[(\Z \Z')^{2}\right]\right)&=&\sum_{k\neq k'}\sum_{i,j}\PE(Z_{i,k}Z_{j,k})\PE(Z_{i,k'}Z_{j,k'})+\sum_{k}\sum_{i}\PE(Z_{i,k}^{2})\\
&=&nN(N-1)+N(N-1)n(n-1)\left(\frac{1}{n-1}\right)^{2}+n^{2}N
\end{eqnarray*}
where the values of the involved expectations may be found in the proof of Lemma \ref{lem:MP_noniid} in Section \ref{tech:lem}. We thus have
$$
\PE\left[\frac{1}{n}\sum_{i=1}^{n}\lambda_{i}^{2}\right]\leq 2 + \frac{n}{N}
$$
which ends the proof.

\subsection{Proof of Theorem \ref{th:consistency}}

The first step is to prove the consistency of $\hat{\eta}$.
We shall first prove that $L_{n}(\eta)$ converges uniformly for $\eta\in[0,1-\delta]$ in probability to $L(\eta)$ given by
$$
L(\eta)=-2\log \sigma^{\star}-\log \int \left[\frac{\eta^{\star}(\lambda-1)+1}{\eta(\lambda-1)+1}\right]  d\mu_{a}(\lambda) - \int \log\left(\eta(\lambda-1)+1\right)d\mu_{a}(\lambda).
$$
Using Lemma \ref{lem:var_quad_form} with $\H$ with diagonal entries $1/(\eta(\lambda_i-1)+1)$, we get that
\begin{eqnarray*}
\Var\left[\frac{1}{n} \sum_{i=1}^{n} \frac{\widetilde{Y}_i^{2}}{\eta(\lambda_i-1)+1}\vert \Z\right]
&\leq &\frac{{\sigma^\star}^4}{n^{2}}\sum_{i=1}^{n}\left[2\left(\frac{\eta^{\star}(\lambda_i-1)+1}{\eta(\lambda_i-1)+1}\right)^{2}+3\left(\frac{1}{q}-1\right)\left(\frac{\eta^{\star}\lambda_i}{\eta(\lambda_i-1)+1}\right)^{2}\right]\\
&\leq &{\sigma^\star}^4\left(2+3\left(\frac{1}{q}-1\right)\right)\frac{1}{n^{2}}\sum_{i=1}^{n}\left(\frac{\lambda_i+1}{\delta}\right)^{2}
\end{eqnarray*}
since $\eta\in[0,1-\delta]$. Now, using Lemma \ref{lem:lambda2} we get that
$$
\frac{1}{n^{2}}\sum_{i=1}^{n}\left(\frac{\lambda_i+1}{\delta}\right)^{2}=o_{P}(1)
$$
which leads to
\begin{eqnarray*}
\frac{1}{n} \sum_{i=1}^{n} \frac{\widetilde{Y}_i^{2}}{\eta(\lambda_i-1)+1}&=&\PE\left[\frac{1}{n} \sum_{i=1}^{n} \frac{\widetilde{Y}_i^{2}}{\eta(\lambda_i-1)+1}\vert \Z\right]+o_{p}(1)\\
&=&{\sigma^{\star}}^{2}\frac{1}{n} \sum_{i=1}^{n} \frac{\eta^{\star}(\lambda_i-1)+1}{\eta(\lambda_i-1)+1}+o_{P}(1).
\end{eqnarray*}
Now, using Lemma \ref{lem:trickMP} we easily get that $\frac{1}{n} \sum_{i=1}^{n} \frac{\eta^{\star}(\lambda_i-1)+1}{\eta(\lambda_i-1)+1}$ converges in probability to $\int [\frac{\eta^{\star}(\lambda-1)+1}{\eta(\lambda-1)+1}]  d\mu_{a}(\lambda)$ and 
$\frac{1}{n}\sum_{i=1}^{n} \log[ (\eta(\lambda_i-1)+1)]$ 
converges in probability to $\int \log (\eta(\lambda-1)+1 )d\mu_{a}(\lambda)$ so that $L_{n}(\eta)=L(\eta)+o_{P}(1)$.

In order to prove the uniform convergence of $L_n$ to $L$ in probability on $[0,1-\delta]$, we shall use the following lemma 
which is proved in section \ref{tech:lem}.
\begin{lemma}\label{lem:cv_unif}
Assume that for any $\eta \in [0, 1 - \delta]$, $L_n(\eta)$ converges in probability to $L(\eta)$ and that 
\begin{equation}
\label{eq:consi1}
\sup_{\eta\in[0,1-\delta]}\left\vert L'_{n}(\eta)\right\vert = O_{P}(1), \textrm{ as $n$ tends to infinity},
\end{equation}
then 
\begin{equation*}
\sup_{\eta\in[0,1-\delta]}\left\vert L_n(\eta)-L(\eta)\right\vert=o_P(1), \textrm{ as $n$ tends to infinity}.
\end{equation*}
\end{lemma}
Let us now prove that $\sup_{\eta\in[0,1-\delta]}\left\vert L'_{n}(\eta)\right\vert = O_{P}(1)$.
Note that
\begin{equation}\label{eq:Ln_prime}
L'_n(\eta)=\left(\frac{1}{n}\sum_{i=1}^n\frac{\widetilde{Y}_i^2(\lambda_i-1)}{\left\{\eta(\lambda_i-1)+1\right\}^2}\right)\left(\frac{1}{n}\sum_{i=1}^n\frac{\widetilde{Y}_i^2}{\eta(\lambda_i-1)+1}\right)^{-1}
-\frac{1}{n}\sum_{i=1}^n\frac{\lambda_i-1}{\eta(\lambda_i-1)+1}.
\end{equation}
A study of $\eta\mapsto \left(\frac{1}{n}\sum_{i=1}^n\frac{\widetilde{Y}_i^2(\lambda_i-1)}{\left\{\eta(\lambda_i-1)+1\right\}^2}\right)\left(\frac{1}{n}\sum_{i=1}^n\frac{\widetilde{Y}_i^2}{\eta(\lambda_i-1)+1}\right)^{-1}
$ shows that it is decreasing and  that it takes negative values for $\eta\in[0,1-\delta]$, so that its absolute value is maximum for $\eta=1-\delta$.
Thus
\begin{multline*}
\sup_{\eta\in[0,1-\delta]}\left\vert L'_{n}(\eta)\right\vert \leq \frac{1}{\delta}
\left(\frac{1}{n}\sum_{i=1}^n \widetilde{Y}_i^2|\lambda_i - 1|\right)\left(\frac{1}{n}\sum_{i=1}^n \widetilde{Y}_i^2\right)^{-1}
+\frac{1}{n\delta}\sum_{i=1}^n |\lambda_i -1|
\\
\leq
\frac{2}{\delta}+\frac{1}{\delta}
\left(\frac{1}{n}\sum_{i=1}^n \widetilde{Y}_i^2\lambda_i \right)\left(\frac{1}{n}\sum_{i=1}^n \widetilde{Y}_i^2\right)^{-1}
+\frac{1}{n\delta}\sum_{i=1}^n \lambda_i
.
\end{multline*}
By Lemma \ref{lem:var_quad_form} with $\H=\Id$, we get
\begin{multline*}
\frac{1}{n}\sum_{i=1}^n \widetilde{Y}_i^2=\PE\left[\frac{1}{n} \sum_{i=1}^{n} \widetilde{Y}_i^{2}\vert \Z\right]+o_{p}(1)
=\frac{\sigma^{2\star}}{n}\sum_{i=1}^n \left[\eta^{\star}(\lambda_i -1)+1)\right] + o_{p}(1) \\
= \sigma^{2\star}\int (\eta(\lambda-1)+1 )d\mu_{a}(\lambda)+ o_{p}(1),
\end{multline*}
where the last equality comes from Lemma \ref{lem:trickMP}. In the same way, we get by using Lemma \ref{lem:var_quad_form} with $\H$ having
its diagonal entries equal to $\lambda_i$ and Lemma \ref{lem:trickMP} that
$$
\frac{1}{n}\sum_{i=1}^n \widetilde{Y}_i^2\lambda_i = \sigma^{2\star}\int \lambda(\eta(\lambda-1)+1 )d\mu_{a}(\lambda)+ o_{p}(1)=O_{P}(1).
$$
Finally, we get from Lemma  \ref{lem:trickMP} that
$$
\frac{1}{n}\sum_{i=1}^n \lambda_i=\int \lambda d\mu_{a}(\lambda)+ o_{p}(1) =O_{P}(1) 
$$
which  ends the proof of \eqref{eq:consi1}. By Lemma \ref{lem:cv_unif}, we thus have proved that
\begin{equation}
\label{eq:consi2}
\sup_{\eta\in[0,1-\delta]}\left\vert L_{n}(\eta)-L(\eta)\right\vert = o_{P}(1).
\end{equation}

Now, using Jensen's inequality, we easily get that for all $\eta\in[0,1]$, $L(\eta)\leq L(\eta^{\star})$, with equality if and only if $\eta=\eta^{\star}$. This together with \eqref{eq:consi2} gives
\begin{equation}
\label{eq:consi3}
\hat{\eta}=\eta^{\star}+o_{P}(1).
\end{equation}

The next step is to prove that $\sqrt{n}(\hat{\eta}-\eta^\star)=O_{P}(1)$. 
Let us first note that $\hat{\eta}$ satisfies the following equation:
\begin{equation}\label{eq:eta_chap_exp}
\sqrt{n}(\hat{\eta}-\eta^\star)=-\frac{\sqrt{n}L'_n(\eta^\star)}{L''_n(\widetilde{\eta})}
,\quad \widetilde{\eta}\in (\hat{\eta},\eta^\star)\;.
\end{equation}
We first prove the asymptotic convergence of $L''_n(\widetilde{\eta})$.

\begin{lemma}\label{lem:LnSeconde}
Let $\Y=(Y_1,\dots,Y_n)'$ satisfy Model (\ref{eq:modele_final}) with $\eta^{\star}>0$ and the entries $W_{i,j}$ of $\W$ satisfy Assumption \ref{assum:A1}. 
Then, for  all $q$ in  $(0,1]$, as $n,N\to\infty$ such that $n/N\to a \in (0,1]$,
for any random variable $\widetilde{\eta}$ such that $\widetilde{\eta}\in (\hat{\eta},\eta^\star)$,
$$
L''_n(\widetilde{\eta})=- {\sigma^{\star}}^{2}\gamma^2(a,\eta^\star)+o_P(1).
$$
\end{lemma}
Lemma \ref{lem:LnSeconde} is proved in Section \ref{tech:lem}.\\
Let us now focus on the properties of $L'_n(\eta^\star)$. Using the following notation
\begin{equation}\label{eq:def:U_i}
U_i=\frac{\widetilde{Y}_i}{\sqrt{\eta^{\star}(\lambda_i-1)+1}}\;,
\end{equation}
we see that $\sqrt{n}L'_n(\eta^{\star})$ can be rewritten as follows:
\begin{align*}
&\left\{\frac{1}{\sqrt{n}}\sum_{i=1}^n\left({U_i}^2-\frac{1}{n}\sum_{j=1}^n U_j^2\right)g(\eta^\star,\lambda_i)\right\}
\left(\frac{1}{n}\sum_{i=1}^n U_i^2\right)^{-1}\\
&=\left\{\frac{1}{\sqrt{n}}\sum_{i=1}^n\left[\left({U_i}^2-1\right)+\left(1-\frac{1}{n}\sum_{j=1}^n U_j^2\right)\right]
g(\eta^\star,\lambda_i)
\right\}\left(\frac{1}{n}\sum_{i=1}^n U_i^2\right)^{-1}\\
&=\left\{\frac{1}{\sqrt{n}}\sum_{i=1}^n \left({U_i}^2-1\right)g(\eta^\star,\lambda_i)\right\}
\left(\frac{1}{n}\sum_{i=1}^n U_i^2\right)^{-1}
-\left\{\frac{1}{\sqrt{n}}\sum_{j=1}^n \left({U_j}^2-1\right)\right\}\left\{\frac{1}{n}\sum_{i=1}^n g(\eta^\star,\lambda_i)\right\}
\left(\frac{1}{n}\sum_{i=1}^n U_i^2\right)^{-1}\;,
\end{align*}
where 
$$
g(\eta,\lambda)=\frac{\lambda-1}{\eta (\lambda-1)+1}\;.
$$
But using Lemma \ref{lem:var_quad_form} and Lemma \ref{lem:trickMP} we get 
$$
\Var\left[n^{-1/2}\sum_{j=1}^n ({U_j}^2-1)\vert \Z\right]=O_P(1)
$$
Moreover, by Lemma \ref{lem:trickMP}, $n^{-1}\sum_{i=1}^n g(\eta^\star,\lambda_i)$ converges in probability
to $\int g(\eta^\star,\lambda)\rmd\mu_a(\lambda)$.
Thus,
\begin{equation}
\label{eq:deriReduit}
\sqrt{n}L'_n(\eta^{\star})=\frac{1}{\sqrt{n}}\sum_{i=1}^n 
\left({U_i}^2-1\right)\left(g(\eta^\star,\lambda_i)-\int g(\eta^\star,\lambda)\rmd\mu_a(\lambda)\right)+o_P(1),
\textrm{ as } n\to\infty\;.
\end{equation}
Using again Lemma \ref{lem:var_quad_form} and Lemma \ref{lem:trickMP} we obtain
$$
\sqrt{n}L'_n(\eta^{\star})=O_{P}(1).
$$
This, together with Lemma \ref{lem:LnSeconde} and \eqref{eq:eta_chap_exp} ends the proof of Theorem \ref{th:consistency}.

\subsection{Proof of Theorem \ref{th:modele_general}}

Notice first that all previous results may be used, replacing  Assumption \ref{assum:A1} by the assumption that the $Z_{i,j}$ are i.i.d. standard Gaussian. Indeed, in this case, Lemma \ref{lem:MP_noniid} reduces to the original result of \cite{marchenko:pastur:1967}, Lemma \ref{lem:trickMP} only involves  Lemma \ref{lem:MP_noniid} and truncation arguments, and the computations leading to Lemma \ref{lem:lambda2} still hold. Thus, Theorem \ref{th:consistency} and Lemma \ref{lem:LnSeconde} also still hold.\\

Let us now prove that $\sqrt{n}L'_n(\eta^{\star})$ converges in distribution to a centered Gaussian.
Define $\H$ the diagonal $n\times n$ matrix with diagonal entries
$$
H_{i}=\frac{1}{\eta^{\star}(\lambda_i-1)+1}\left[g(\eta^\star,\lambda_i)-\int g(\eta^\star,\lambda)\rmd\mu_a(\lambda)\right].
$$
Define
$$
\L_{n}=\frac{1}{\sqrt{n}}\widetilde{\Y}'\H\widetilde{\Y}.
$$
Then using  \eqref{eq:deriReduit} and Lemma \ref{lem:trickMP} we have
$$
\sqrt{n}L'_n(\eta^{\star})=\L_{n}-\PE[\L_{n}\vert \Z]+o_P(1).
$$
Now using Lemma \ref{lem:var_quad_form} we get that
setting
$\gamma_n^{2}=\Var\left[\L_{n}\vert  \Z\right]$,
\begin{multline*}
\gamma_n^{2}=2{\sigma^{\star}}^{4}\frac{1}{n}\Tr\left[\H^{2}\left((1-\eta)^{\star}Id_{\rset^{n}}+\eta^{\star}\D\right)^{2}\right]+3{\sigma^{\star}}^{4}{\eta^{\star}}^{2}\left(\frac{1}{q}-1\right)\frac{1}{n}\sum_{i=1}^{N}M_{i,i}^{2}\\
=2{\sigma^{\star}}^{4}\frac{1}{n}\sum_{i=1}^{n}\left(g(\eta^\star,\lambda_i)-\int g(\eta^\star,\lambda)\rmd\mu_a(\lambda)\right)^{2}\\
+3{\sigma^{\star}}^{4}{\eta^{\star}}^{2}\left(\frac{1}{q}-1\right)\frac{1}{n}\sum_{i=1}^{n}\sum_{k,l=1}^{n}\lambda_{k}\lambda_{l}H_{k}H_{l}V_{i,k}^{2}V_{i,l}^{2}.
\end{multline*}
The first term in this sum converges as $n,N\to \infty$ to $2{\sigma^{\star}}^{4}\gamma^2(a,\eta^\star)$.\\
Under the assumption that the $Z_{i,j}$ are i.i.d. standard Gaussian, the matrix of eigenvectors $\V$ is Haar distributed on the orthonormal matrices, and is independent of 
$(\lambda_{i})_{1\leq i \leq n}$, see
\cite{Bai:Silverstein:book} chapter 10. Conditionally to the eigenvalues $(\lambda_{i})_{1\leq i \leq n}$, we thus get that
$$
\PE\left[\frac{1}{n}\sum_{i=1}^{n}\sum_{k,l=1}^{n}\lambda_{k}\lambda_{l}H_{k}H_{l}V_{i,k}^{2}V_{i,l}^{2}\vert\D\right]=\left(\frac{1}{N}\sum_{k=1}^{n}\lambda_{k}H_{k}\right)^{2}(1+o(1))
$$
converges to 
$$
a^{2}\left[\int \frac{\lambda(\lambda-1)}{(\eta^\star(\lambda-1)+1)^{2}}\rmd\mu_a(\lambda)-\int\frac{\lambda}{(\eta^\star(\lambda-1)+1)}\rmd\mu_a(\lambda)\int\frac{\lambda-1}{(\eta^\star(\lambda-1)+1)}\rmd\mu_a(\lambda)\right]^{2}
$$
and 
$$
\Var\left[\frac{1}{n}\sum_{i=1}^{n}\sum_{k,l=1}^{n}\lambda_{k}\lambda_{l}H_{k}H_{l}V_{i,k}^{2}V_{i,l}^{2}\vert \D\right]=o_{P}(1)
$$
so that 
$$\gamma_n^{2}=2{\sigma^{\star}}^{4}\gamma^2(a,\eta^\star)+3{\sigma^{\star}}^{4}{\eta^{\star}}^{2}\left(\frac{1}{q}-1\right)S(a,\eta^\star)
+o_{P}(1).
$$
Denote $\Delta$ the diagonal $N\times N$-matrix with diagonal entries $\Delta_{i}=\frac{\sigma^{\star}\sqrt{\eta^{\star}}}{\sqrt{q}}\pi_{i}$.
Let us now write
$$
\L_{n}-\PE(\L_{n}\vert \Z)=\L_{n}-\PE\left[\L_{n}\vert \Delta, \Z\right]+\PE\left[\L_{n}\vert \Delta, \Z\right] -\PE\left[\L_{n}\vert  \Z\right].
$$
We first have
$$
\PE\left[\L_{n}\vert \Delta, \Z\right] -E\left[\L_{n}\vert  \Z\right]={\sigma^{\star}}^{2}\eta^{\star}\frac{1}{\sqrt{n}}\sum_{i=1}^{N}\left(\frac{\pi_{i}^{2}}{q}-1\right)M_{i,i}
$$
whose variance, conditionally to $\Z$ is
$$
s_{n,1}^{2}={\sigma^{\star}}^{4}{\eta^{\star}}^{2}\left(\frac{1}{q}-1\right)\frac{1}{n}\sum_{i=1}^{N}M_{i,i}^{2}.
$$
In the same way as for $\gamma_n^{2}$ we get that
$$
s_{n,1}^{2}={\sigma^{\star}}^{4}{\eta^{\star}}^{2}\left(\frac{1}{q}-1\right)S(a,\eta^\star)
+o_{P}(1).
$$
Let
$$
\xi_{i}=\left(\frac{\pi_{i}^{2}}{q}-1\right)M_{i,i}=\left(\frac{\pi_{i}^{2}}{q}-1\right)\sum_{k=1}^{n}\frac{\lambda_{k}(\lambda_{k}-1)}{(\eta^\star(\lambda_{k}-1)+1)^{2}}V_{i,k}^{2}.
$$
Since $\eta^{\star} >0$, the function $\lambda \mapsto \frac{\lambda(\lambda-1)}{(\eta^\star(\lambda-1)+1)^{2}}$ is bounded, and $\sum_{k=1}^{n}V_{i,k}^{2}\leq 
\sum_{k=1}^{N}V_{i,k}^{2}=1$. Also, the variables $\left(\frac{\pi_{i}^{2}}{q}-1\right)$ are uniformly bounded by $1/q$. Thus
$$
\frac{1}{n}\sum_{i=1}^{n} 
\PE\left[\xi_{i}^{2}\1_{|\xi_{i}|\geq cn} \vert \Z\right]=0
$$
for large enough $n$.
Then, by Lindeberg's Theorem, conditionally to $\Z$,
$$
\frac{1}{s_{n,1}}\left(\PE\left[\L_{n}\vert \Delta, \Z\right] -\PE\left[\L_{n}\vert  \Z\right]\right)
$$
converges in distribution to $\mathcal{N}(0,1)$.
\\
Let us now set
$$
s_{n,2}^{2}=\gamma_n^{2}-s_{n,1}^{2}
$$
and notice that $s_{n,2}^{2}$ converges to 
$$
2{\sigma^{\star}}^{4}\gamma^2(a,\eta^\star)+2{\sigma^{\star}}^{4}{\eta^{\star}}^{2}\left(\frac{1}{q}-1\right)S(a,\eta^\star).
$$
We shall prove that, conditionally to $\Z$ and $\Delta$, $(\L_{n}-\PE(\L_{n}\vert \Delta, \Z))/s_{n,2}$ converges in distribution to $\mathcal{N}(0,1)$, and thus also
unconditionally. Write
$$
\L_{n}=\frac{1}{\sqrt{n}}\begin{array}{ll}(\w'&\eps')\\
&\end{array}
B\left(\begin{array}{l}\w\\
\eps
\end{array}\right)
$$
where $B$ is the $(N+n)\times (N+n)$-matrix 
$$
B=\left(\begin{array}{cc}\Delta & 0\\
0 & \sigma^{\star}\sqrt{(1-\eta)^{\star}} Id_{\rset^{n}}
\end{array}
\right)
\left(\begin{array}{cc}
\V\left(\begin{array}{cc}\D\H & 0\\
0 & 0
\end{array}\right)\V' \;&\;\tilde{\V}\sqrt{\D} \H\\
\begin{array}{l}\\ \H\sqrt{\D}\tilde{\V}'\end{array}& \begin{array}{l}\\ \H\end{array}
\end{array}\right) 
\left(\begin{array}{cc}\Delta & 0\\
0 & \sigma^{\star}\sqrt{(1-\eta)^{\star}} Id_{\rset^{n}}
\end{array}
\right).
$$
Here, $\tilde{\V}$ is the $N\times n$-matrix which consists of the first $n$ columns of $\V$.
Let $\phi$ be the characteristic function of
$(\L_{n}-\PE(\L_{n}\vert \Delta, \Z))/s_{n,2}$ conditionally to $\Z$ and $\Delta$. Notice first that if $b_{j}$, $j=1,\ldots,n+N$ are the eigenvalues of $B$, we may write
$$
\L_{n}-\PE\left[\L_{n}\vert \Delta, \Z\right]=\frac{1}{\sqrt{n}}\sum_{j=1}^{N+n}b_{j}(e_{j}^{2}-1).
$$
for random variables $e_{j}$ i.i.d. standard Gaussian.
Thus
$$
\phi \left(t\right)=\prod_{j=1}^{N+n}\left[\left(1-2i\frac{tb_{j}}{s_{n,2}\sqrt{n}}\right)^{-1/2}\exp\left(-i\frac{tb_{j}}{s_{n,2}\sqrt{n}}\right)\right]
$$
and Taylor expansion leads to 
\begin{eqnarray*}
\log \phi \left(t\right)&=&\sum_{j=1}^{N+n}\left[-\frac{1}{2}\log \left(1-2i\frac{tb_{j}}{s_{n,2}\sqrt{n}}\right)-i\frac{tb_{j}}{s_{n,2}\sqrt{n}}\right]\\
&=&-t^{2}\frac{1}{ns_{n,2}^{2}}\sum_{j=1}^{N+n}b_{j}^{2}+O\left[\frac{1}{n\sqrt{n}s_{n,2}^{3}}\sum_{j=1}^{N+n}b_{j}^{3}\right].
\end{eqnarray*}
We shall now prove that $\frac{1}{ns_{n,2}^{2}}\sum_{j=1}^{N+n}b_{j}^{2}$ converges to $1/2$.
Tedious computations give
\begin{eqnarray*}
\sum_{j=1}^{N+n}b_{j}^{2}&=&\Tr(B^{2})\\
&=&\Tr(\Delta M\Delta^{2}M\Delta)+{\sigma^{\star}}^{4}(1-{\eta^{\star}}^{2}\Tr (\H^{2})+2{\sigma^{\star}}^{2}(1-{\eta^{\star}})\Tr[\Delta^{2}\tilde{\V}\D\H^{2}\tilde{\V}'].
\end{eqnarray*}
Using the distribution of $\V$ and its independence on $\D$ we get
\begin{eqnarray*}
\PE\left[\sum_{j=1}^{N+n}b_{j}^{2}\vert \D\right]&=&2{\sigma^{\star}}^{4}\Tr\left[\H^{2}\left((1-\eta)^{\star}Id_{\rset^{n}}+\eta^{\star}\D\right)^{2}\right]\\
&&+2{\sigma^{\star}}^{4}{\eta^{\star}}^{2}\left(\frac{1}{q}-1\right)\left(\frac{1}{N}\sum_{k=1}^{n}\lambda_{k}H_{k}\right)^{2}(1+o(1))
\end{eqnarray*}
so that 
$$\PE\left[\frac{1}{n}\sum_{j=1}^{N+n}b_{j}^{2}\vert \D\right]=2{\sigma^{\star}}^{4}\gamma^2(a,\eta^\star)+2{\sigma^{\star}}^{4}{\eta^{\star}}^{2}\left(\frac{1}{q}-1\right)S(a,\eta^\star)
+o_{P}(1).
$$ 
Moreover, tedious computations again give
$$
\Var\left[\frac{1}{n}\sum_{j=1}^{N+n}b_{j}^{2}\vert \D\right]=o_{P}(1),
$$
and we obtain that 
$$
\frac{1}{ns_{n,2}^{2}}\sum_{j=1}^{N+n}b_{j}^{2}=\frac{1}{2}+o_{P}(1).
$$
\\
We shall now prove that $\frac{1}{n\sqrt{n}s_{n,2}^{3}}\sum_{j=1}^{N+n}b_{j}^{3}=o_{P}(1)$. To do so, it is enough to prove that
$\max_{j}|b_{j}|=o_{P}(\sqrt{n})$. Notice that for any normed vector $A=(A_{1},A_{2})$ in $\rset^{N+n}$ where $A_{1}\in\rset^{N}$ and $A_{2}\in\rset^{n}$,
$$\max_{j}|b_{j}|\leq A'BA.
$$
Now,
$$
A'BA=A_{1}'(\Delta M \Delta)A_{1}+2{\sigma^{\star}}\sqrt{1-\eta^{\star}} A_{1}'(\Delta\tilde{\V} \sqrt{\D}\H)A_{2} + {\sigma^{\star}}^{2}(1-\eta^{\star})A_{2}'\H A_{2}.
$$
First, since $\eta^{\star}>0$, all entries of $\H$ and $\D$ and $\H\D$ are uniformly bounded and so are all entries of $\Delta$. We thus get $A_{2}'\H A_{2}=O(1)$
and $A_{1}'(\Delta\tilde{\V} \sqrt{\D}\H)A_{2}=O(1)$. Then, using the distribution of $\V$ and its independence on $\D$ we get
$$\PE\left[A_{1}'(\Delta M \Delta)A_{1}\vert \D\right]=O\left(\frac{1}{N}\sum_{i=1}^{n}\lambda_{i}H_{i}
\right)
$$
and
$$\Var\left[A_{1}'(\Delta M \Delta)A_{1}\vert \D\right] = o_{P}(1),
$$
so that $A'BA=O_{P}(1)$. We have thus proved that $\max_{j}|b_{j}|=O_{P}(1)=o_{P}(\sqrt{n})$. 
\\

Thus $\phi(t)$ converges in probability for all $t$ to $\exp - \frac{t^{2}}{2}$ and the convergence may be strengthened by contradiction to an a.s. convergence,  so that  conditionally to $\Z$ and $\Delta$, $(\L_{n}-\PE(\L_{n}\vert \Delta, \Z))/s_{n,2}$ converges in distribution to 
$\mathcal{N}(0,1)$.\\

Now, conditionally to $\Z$ and $\Delta$, $(\L_{n}-E(\L_{n}\vert \Delta, Z))/s_{n,2}$ converges in distribution to a Gaussian random variable independent of $\Delta$. Thus
conditionally to $\Z$, $\L_{n}-E\left[\L_{n}\vert \Delta, Z\right]$ and $E\left[\L_{n}\vert \Delta, Z\right] -E\left[\L_{n}\vert  Z\right]$ converge in distribution to independent Gaussian variables, so that their sum converges in distribution to a centered Gaussian with variance the sum of the variances, namely the limit of $\gamma_n^{2}$, and Theorem
\ref{th:modele_general} is proved.

\subsection{Proof of Theorem \ref{th:modelenonsparse}} 
Using Lemma \ref{lem:LnSeconde} and \eqref{eq:eta_chap_exp} , there remains to prove that $\sqrt{n}L'_n(\eta^{\star})$ converges in distribution to 
$\mathcal{N}(0,2{\sigma^{\star}}^{4}\gamma^2(a,\eta^\star))$ and that $\gamma_{n}^{2}$ converges in probability to $\gamma^2(a,\eta^\star)$.\\
Notice first that when $q=1$, $(U_{1},\ldots,U_{n})\vert \Z$ is a centered Gaussian vector with a covariance matrix equal to ${\sigma^{\star}}^{2}$ times the identity matrix. We shall prove that conditionally to $\Z$, $\sqrt{n}L'_n(\eta^{\star})$ converges in distribution to 
$\mathcal{N}(0,2{\sigma^{\star}}^{4}\gamma^2(a,\eta^\star))$ so that the result still holds unconditionally. Using \eqref{eq:deriReduit}, it is only needed to prove it for
$\frac{1}{\sqrt{n}}\sum_{i=1}^n 
\left({U_i}^2-1\right)\left(g(\eta^\star,\lambda_i)-\int g(\eta^\star,\lambda)\rmd\mu_a(\lambda)\right)$. Now, conditionally to $\Z$, the variance of
$$
\sum_{i=1}^n 
\left({U_i}^2-1\right)\left(g(\eta^\star,\lambda_i)-\int g(\eta^\star,\lambda)\rmd\mu_a(\lambda)\right)
$$
is 
$$
\gamma_n^{2}=\frac{2{\sigma^{\star}}^{4}}{n}\sum_{i=1}^{n}\left(g(\eta^\star,\lambda_i)-\int g(\eta^\star,\lambda)\rmd\mu_a(\lambda)\right)^{2}.
$$
Since $\eta^{\star} >0$, $g(\eta^\star,\lambda)$ is a bounded function of $\lambda$, and using Lemma \ref{lem:trickMP}, 
$$
\gamma_n^{2}=2{\sigma^{\star}}^{4}\gamma^2(a,\eta^\star))+o_{P}(1).
$$
Also, setting $\xi_{i}=\left({U_i}^2-1\right)\left(g(\eta^\star,\lambda_i)-\int g(\eta^\star,\lambda)\rmd\mu_a(\lambda)\right)$ and $C$ an upper bound of $|g(\eta^\star,\lambda)|$, we get that for any $c>0$,
\begin{eqnarray*}
\frac{1}{n}\sum_{i=1}^{n} 
\PE\left[\xi_{i}^{2}\1_{|\xi_{i}|\geq cn} \vert \Z\right]&\leq&4C^{2}{\sigma^{\star}}^{4}\PE \left[\left({U_1}^2-1\right)^{2}\1_{2C|{U_1}^2-1|\geq cn} \vert Z\right] \\
&=& 4C^{2}{\sigma^{\star}}^{4}\PE \left[\left({U_1}^2-1\right)^{2}\1_{2C|{U_1}^2-1|\geq cn} \right]
=o(1),
\end{eqnarray*}
where the first equality comes from the fact that the distribution of $(U_{1},\ldots,U_{n})\vert \Z$  does not depend on $\Z$
and is thus also the distribution of $(U_{1},\ldots,U_{n})$.
Then,
using Lindeberg's Theorem, conditionally to $\Z$, $\sqrt{n}L'_n(\eta^{\star})$ converges in distribution to 
$\mathcal{N}(0,2{\sigma^{\star}}^{4}\gamma^2(a,\eta^\star))$ and thus also unconditionally.\\
The fact that $\gamma_{n}^{2}$ converges in probability to $\gamma^2(a,\eta^\star)$ is a straightforward consequence of Taylor expansion, the fact that 
 $g(\eta^\star,\lambda)$ and its derivative with respect to $\eta$ in the neighborhood of $\eta^{\star}$ are bounded functions of $\lambda$, and Slutzky's Lemma.  

\subsection{Proofs of technical lemmas}\label{tech:lem}

\subsubsection{Proof of Lemma \ref{lem:MP_noniid}}
As a byproduct of Theorem 1.1, Corollary 1.1 and Remark 1.1 of \cite{bai:zhou:2008}, we use the following result to prove Lemma  \ref{lem:MP_noniid}.
 \begin{theorem2}[Bai and Zhou (2008)]
Let $\Z$ be a matrix of size $n \times N$ which columns, denoted by $Z_1,\dots,Z_N$, are independent and let us denote $\bar{\Z}=\frac{1}{N}\sum_{k=1}^{N} Z_k$. Let us also recall that $\R=\Z\Z'/N$ and $F^{\R}$ is its empirical spectral distribution defined by $F^{\R}(x)=\frac{1}{n}\sum_{k=1}^{n} 1_{\{\lambda_k\geq x\}}$, where $\lambda_1,\dots,\lambda_n$ are the eigenvalues of $\R$.
As $N \to \infty$, assume the following:
\begin{enumerate}
\item $T=(t_{i,j})$ is a matrix such that $\PE(\bar{Z}_{i,j} Z_{m,j})=t_{m,i}$ for all $j$ .
\item $\frac{1}{N}\underset{i \neq m}{\max}$ $\PE(\bar{Z}_{i,j} Z_{m,j})^2 \to 0$ uniformly in $j \leq N$.
\item $\frac{1}{N^2}\sum_{\Lambda}\left(\PE(\bar{Z}_{i,j} Z_{m,j}-t_{m,i})(Z_{i',j} \bar{Z}_{m',j}-t_{i',m'})\right)^2 \to 0$ uniformly in $j \leq N$, with
 $\Lambda=\{(i,m,i',m') : 1 \leq i,m,i',m' \leq n \} \backslash \{(i,m,i',m') : i=i' \neq m=m' $ or $ i=m' \neq i'=m \}$.
\item $\frac{n}{N} \to a \in (0,+\infty)$.
\item The norm of $T$ is uniformly bounded and $F^T$ tends to a degenerate distribution with mass at $1/a$.
\end{enumerate}
Then, with probability 1, $F^{R}$ converges to the Marchenko-Pastur distribution defined in (\ref{eq:marchenko}).
\end{theorem2}

Observe that for all $j=1,\dots,N$,
\begin{equation}\label{eq:sum_Z}
\sum_{i=1}^{n}Z_{i,j}=0
\end{equation}
and
\begin{equation}\label{eq:sum_Z_car}
\sum_{i=1}^{n}Z_{i,j}^{2}=n.
\end{equation}
Moreover, for each $j$, the random variables $(Z_{i,j})_{1\leq i\leq n}$ are exchangeable. Thus, we deduce from
(\ref{eq:sum_Z_car}) that for all $i=1,\dots,n$ and $j=1,\dots,N$,
$
\PE(Z_{i,j}^2)=1.
$
Hence, by (\ref{eq:sum_Z}), we get that
$$
0=\left(\sum_{i=1}^{n}Z_{i,j}\right)^2=\sum_{i=1}^n Z_{i,j}^2+\sum_{1\leq i\neq m\leq n} Z_{i,j} Z_{m,j}\;,
$$
which, by (\ref{eq:sum_Z_car}), implies that for all $j=1,\dots,N$ and $i\neq m=1,\dots,n$,
\begin{equation}\label{eq;prod_Z}
\PE(Z_{i,j} Z_{m,j})=-\frac{n}{n(n-1)}=-\frac{1}{n-1}\;.
\end{equation}
Thus, the matrix $\T=\T_n$ defined in Theorem (Bai and Zhou (2008)) is equal to
$
\T=n/(n-1)\Id_{\rset^n}-\J_n/(n-1)\;,
$ 
where $\J_n$ is a $n\times n$ matrix having all its entries
equal to 1. Hence the eigenvalues of $\T$ are 0 with multiplicity 1 and $n/(n-1)$ with multiplicity
$(n-1)$, which gives Condition 5. of Theorem (Bai and Zhou (2008)). 

Let us then check Condition 2. of Theorem (Bai and Zhou (2008)). 
Observe that, for $i\neq m$, $\PE[(Z_{i,j}Z_{m,j}-t_{m,i})^2]=\PE(Z_{i,j}^2Z_{m,j}^2)-t_{m,i}^2$.
By (\ref{eq:sum_Z_car}), for all $j=1,\dots,N$, 
$$
n^2=\left(\sum_{i=1}^{n}Z_{i,j}^{2}\right)^2=\sum_{i=1}^n Z_{i,j}^{4}+\sum_{1\leq i\neq m\leq n} Z_{i,j}^2 Z_{m,j}^2\;.
$$
Since the $(Z_{i,j})_{1\leq i\leq n}$ are exchangeable for each $j=1,\dots,N$, we get that for all $j=1,\dots,N$,
$$
n=\PE[Z_{1,j}^4]+(n-1)\PE[Z_{1,j}^2 Z_{2,j}^2]\;.
$$
Thus,  for all $j=1,\dots,N$, $\PE[Z_{1,j}^2 Z_{2,j}^2]\leq n/(n-1)$, which with the definition of the $t_{m,i}$'s gives the 
result.

Let us now check Condition 3. of Theorem (Bai and Zhou (2008)). 
Since the random variables $(Z_{i,j})_{1\leq i\leq n}$ are exchangeable, it is enough to prove that, uniformly in $k$,
\begin{enumerate}[(i)]
\item $\PE[Z_{1,k}^{4}]=o(\sqrt{n})$, 
\item $\PE[Z_{1,k}^{2}Z_{2,k}^{2}]-1=o(1)$,
\item $\PE[Z_{1,k}^{3}Z_{2,k}]=o(1)$,
\item $\sqrt{n}\PE[Z_{1,k}^{2}Z_{2,k}Z_{3,k}]=o(1)$,
\item $n\PE[Z_{1,k}Z_{2,k}Z_{3,k}Z_{4,k}]=o(1)$ \textrm{, as $n\to\infty$.}
\end{enumerate}
Observe that (i) implies (ii). Using (\ref{eq:sum_Z}), by expanding
$0=\left(\sum_{i=1}^{n}Z_{i,k}\right)^2 \left(\sum_{i=1}^{n}Z_{i,k}^{2}\right)$ and taking the expectation, 
we get that (i) and (iii) imply (iv). By expanding $0=\left(\sum_{i=1}^{n}Z_{i,k}\right)^4$, which comes
from (\ref{eq:sum_Z}), and by taking the expectation, (i) and (iii) imply (v).
Hence, it is enough to prove (i) and (iii) to conclude the proof of Lemma \ref{lem:MP_noniid}.

Let us first prove (i). By the definition of $Z_{1,k}$ given in (\ref{eq:normalization_1}), we get that for all $k$,
$Z_{1,k}^2\leq n$. Hence,
$$
Z_{1,k}^2\leq \frac{(W_{1,k}-\overline{W}_k)^2}{2\sigma_k^2}\1_{\{s_k^2\geq 2\sigma_k^2\}}+ n\1_{\{s_k^2> 2\sigma_k^2\}}\;,
$$
and, by the assumptions on the $W_{i,k}$'s and on the $\sigma_k$'s,
$$
\PE(Z_{1,k}^4)\leq \frac{W_M^2}{2\kappa^2}+2n^2\PP(s_k^2-\sigma_k^2> \sigma_k^2)\;.
$$
Theorem A of \cite[p. 201]{serfling:1980} implies that the second term of the previous
inequality tends to zero as $n$ tends to infinity uniformly in $k$, which concludes the proof of (i).

Let us now prove (iii). Using (\ref{eq:sum_Z}), we get
$
Z_{1,k}^3\left(\sum_{i=1}^n Z_{i,k}\right)=0.
$
By expanding this equation and taking the expectation, we obtain that
$
\PE(Z_{1,k}^4)+\sum_{i=2}^n \PE(Z_{1,k}^3 Z_{i,k})=0.
$
Since the $(Z_{i,k})_{1\leq i\leq n}$ are exchangeable: $\PE(Z_{1,k}^3 Z_{2,k})=-\PE(Z_{1,k}^4)/(n-1)=o(n^{-1/2})$, where the
last equality comes from (i).


\subsubsection{Proof of Lemma \ref{lem:var_quad_form}}
Using (\ref{eq:modele_svd}) and the independence assumptions, we get
\begin{multline}\label{eq:var_form_quad}
\Var(\widetilde{\Y}'\H\widetilde{\Y}|\Z)=\Var\left[\v'\V
\left(
\begin{matrix}
\D\H & \mathbf{0}\\
\mathbf{0} & \mathbf{0}
\end{matrix}
\right)
\V'\v+2\sigma^\star\sqrt{1-\eta^\star}\v'\V\left(
\begin{matrix}
\sqrt{\D}\\
\mathbf{0}
\end{matrix}
\right)
\H\widetilde{\eps}
+{\sigma^\star}^2 (1-\eta^\star)\widetilde{\eps}'\H\widetilde{\eps}
|\Z\right]\\
=\Var\left[\v'\M\v|\Z\right]+4{\sigma^\star}^2 (1-\eta^\star)\Var\left[\v'\V\left(
\begin{matrix}
\sqrt{\D}\\
\mathbf{0}
\end{matrix}
\right)
\H\widetilde{\eps}|\Z\right]+2{\sigma^\star}^4(1-\eta^\star)^2\Tr(\H^2)\;,
\end{multline}
where $\M=\V\left(
\begin{matrix}
\D\H & \mathbf{0}\\
\mathbf{0} & \mathbf{0}
\end{matrix}
\right)\V'$. Using the independence assumptions, we get that
\begin{equation}\label{eq:term2}
4{\sigma^\star}^2 (1-\eta^\star)\Var\left[\v'\V\left(
\begin{matrix}
\sqrt{\D}\\
\mathbf{0}
\end{matrix}
\right)\H\widetilde{\eps}|\Z\right]
=4{\sigma^\star}^4 \eta^\star(1-\eta^\star)\Tr(\B\B')=4{\sigma^\star}^4 \eta^\star(1-\eta^\star)\Tr(\D\H^2)\;,
\end{equation}
where $\B=\V\left(
\begin{matrix}
\sqrt{\D}\\
\mathbf{0}
\end{matrix}
\right)\H$. Moreover, $\PE(\v'\M\v|\Z)={\sigma^\star}^2\eta^\star\Tr(\D^2\H^2)$ and 
\begin{multline*}
\PE\left[(\v'\M\v)^2|\Z\right]=\frac{{\sigma^\star}^4{\eta^\star}^2}{q^2}\left[2 q^2\sum_{1\leq i\neq j\leq N} M_{ij}^2
+q^2\sum_{1\leq i\neq i'\leq N} M_{ii}M_{i'i'}+3q\sum_{1\leq i\leq N} M_{ii}^2\right]\\
={\sigma^\star}^4{\eta^\star}^2\left[2\Tr(\M^2)-2\sum_{1\leq i\leq N} M_{ii}^2 +\Tr(\M)^2-\sum_{1\leq i\leq N} M_{ii}^2
+\frac{3}{q}\sum_{1\leq i\leq N} M_{ii}^2\right]\\
={\sigma^\star}^4{\eta^\star}^2\left[2\Tr(\D^2\H^2) +\Tr(\M)^2+3\left(\frac{1}{q}-1\right)\sum_{1\leq i\leq N} M_{ii}^2\right]\;.
\end{multline*}
Thus,
\begin{equation}\label{eq:term1}
\Var\left[\v'\M\v\right|\Z]={\sigma^\star}^4{\eta^\star}^2\left[2\Tr(\D^2\H^2)+3\left(\frac{1}{q}-1\right)\sum_{1\leq i\leq N} M_{ii}^2\right]\;.
\end{equation}
The proof of the equality in Lemma \ref{lem:var_quad_form} follows from (\ref{eq:var_form_quad}), (\ref{eq:term2}) and (\ref{eq:term1}).
The proof of the inequality in Lemma \ref{lem:var_quad_form} follows now from 
$$
\sum_{1\leq i\leq N} M_{ii}^2\leq \sum_{1\leq i,j\leq N} M_{ij}^2 = \Tr[\D^{2}\H^{2}].
$$
 \subsubsection{Proof of Lemma \ref{lem:cv_unif}}
Let $\epsilon >0$ and let \{$\eta_1 < \dots < \eta_{K(\epsilon)}$\} be a grid of $[0, 1-\delta]$ such that $\vert \eta_j - \eta_{j+1} \vert < \epsilon$  for all $j \in \{0,\dots,K_{\epsilon}\}$ then
\begin{align*}
\sup_{\eta \in [0, 1-\delta]} \vert L_n(\eta) - L(\eta) \vert
&\leq \sup_{j \in \{1, \dots, K_{\epsilon}\}} \left[\sup_{\eta' \in [\eta_j,\eta_{j+1}]} \vert L_n(\eta') - L_n(\eta_j) \vert 
+ \vert L_n(\eta_j) - L(\eta_j) \vert\right.\\
& \left. + \sup_{\eta' \in [\eta_j,\eta_{j+1}]} \vert L(\eta_j) - L(\eta')\vert\right]\\
&\leq \epsilon \sup_{\eta \in [0, 1-\delta]} \vert L_n'(\eta) \vert
+\sup_{j \in \{1, \dots, K_{\epsilon}\}}  \vert L_n(\eta_j) - L(\eta_j) \vert + \omega(\epsilon),
\end{align*}
where $\omega(\epsilon)$ is the modulus of continuity of $L$, which is continuous on the compact $[0, 1 -\delta]$ 
and thus uniformly continuous on this compact.
Since $\underset{\eta \in [0, 1-\delta]}{\sup}\left\vert L'_{n}(\eta)\right\vert = O_{P}(1)$ then, for every $\beta>0$, there exists $C$ 
such that for all $n$, $\PP(\underset{\eta \in [0, 1-\delta]}{\sup}\vert L_n'(\eta) \vert \geq C)\leq\beta$.
Let $\alpha >0$ and let us consider the $\epsilon$-grid such that $\epsilon \leq \alpha/3C$ and $\omega(\epsilon) \leq \alpha/3$, thus we get
that
\begin{align*}
&\PP(\underset{\eta \in [0, 1-\delta]}{\sup}\left\vert L_n(\eta) - L(\eta) \right\vert \geq \alpha) \\
&\leq \PP(\underset{\eta \in [0, 1-\delta]}{\sup} \left\vert L_n'(\eta) \right\vert  \geq C) 
 + \PP( \sup_{j \in \{1, \dots, K_{\epsilon}\}} \left\vert L_n(\eta_j) - L(\eta_j) \right\vert\geq \alpha - C \epsilon - \omega(\epsilon)) \\
& \leq \PP(\underset{\eta \in [0, 1-\delta]}{\sup} \left\vert L_n'(\eta) \right\vert  \geq C) 
+ \PP(\sup_{j \in \{1, \dots, K_{\epsilon}\}} \left\vert L_n(\eta_j) - L(\eta_j) \right\vert \geq \frac{\alpha}{3})\\
& \leq \PP(\underset{\eta \in [0, 1-\delta]}{\sup} \left\vert L_n'(\eta) \right\vert  \geq C) 
+\sum_{j=1}^{K_\epsilon}  \PP( \left\vert L_n(\eta_j) - L(\eta_j) \right\vert \geq \frac{\alpha}{3}),
\end{align*}
which concludes the proof of Lemma \ref{lem:cv_unif} since each term tends to zero as $n$ tends to infinity.

\subsubsection{Proof of Lemma \ref{lem:LnSeconde}}
The second derivative of $L_{n}$ is given by
\begin{align}\label{eq:Ln_second}
L''_n(\eta)&=\left(-\frac2n\sum_{i=1}^n\frac{\widetilde{Y}_i^2(\lambda_i-1)^2}{\left\{\eta(\lambda_i-1)+1\right\}^3}\right)
\left(\frac1n\sum_{i=1}^n\frac{\widetilde{Y}_i^2}{\left\{\eta(\lambda_i-1)+1\right\}}\right)^{-1}\nonumber\\
&+\left(\frac1n\sum_{i=1}^n\frac{\widetilde{Y}_i^2(\lambda_i-1)}{\left\{\eta(\lambda_i-1)+1\right\}^2}\right)^2
\left(\frac1n\sum_{i=1}^n\frac{\widetilde{Y}_i^2}{\left\{\eta(\lambda_i-1)+1\right\}}\right)^{-2}
+\frac1n\sum_{i=1}^n\frac{(\lambda_i-1)^2}{\left\{\eta(\lambda_i-1)+1\right\}^2}\;.
\end{align}
In particular for $\eta=\eta^{\star}$, we have
$$
\frac1n\sum_{i=1}^n\frac{\widetilde{Y}_i^2}{\left\{\eta^{\star}(\lambda_i-1)+1\right\}}=1+o_{P}(1),
$$
and using as previously Lemma \ref{lem:var_quad_form}, Lemma \ref{lem:trickMP} and the fact that all functions of $\lambda$ involved in the empirical means are bounded since $\eta^{\star}>0$, we get
\begin{eqnarray*}
\frac2n\sum_{i=1}^n\frac{\widetilde{Y}_i^2(\lambda_i-1)^2}{\left\{\eta(\lambda_i-1)+1\right\}^3}&=&
\frac{{2\sigma^{\star}}^{2}}{n}\sum_{i=1}^n\frac{(\lambda_i-1)^2}{\left\{\eta(\lambda_i-1)+1\right\}^2} +o_{P}(1)\\
&=&{2\sigma^{\star}}^{2}\int \frac{(\lambda-1)^2}{\left\{\eta(\lambda-1)+1\right\}^2} d\mu_{a}(\lambda) +o_{P}(1)
\end{eqnarray*}
and
\begin{eqnarray*}
\frac1n\sum_{i=1}^n\frac{\widetilde{Y}_i^2(\lambda_i-1)}{\left\{\eta(\lambda_i-1)+1\right\}^2}&=&
\frac{{\sigma^{\star}}^{2}}{n}\sum_{i=1}^n\frac{(\lambda_i-1)}{\left\{\eta(\lambda_i-1)+1\right\}} +o_{P}(1)\\
&=&{\sigma^{\star}}^{2}\int \frac{(\lambda-1)}{\left\{\eta(\lambda-1)+1\right\}} d\mu_{a}(\lambda) +o_{P}(1)
\end{eqnarray*}
leading to
$$
L''_n(\eta)=-{\sigma^{\star}}^{2}\gamma^2(a,\eta^\star)+o_{P}(1).
$$
Using Slutzky's Lemma and $\hat{\eta}=\eta^{\star}+o_{P}(1)$, there just remains to prove that for small enough $\alpha >0$,
$$
\sup_{|\eta-\eta^{\star}|\leq \alpha} |L''_n(\eta)-L''_n(\eta)|=O_{p}(\alpha).
$$
But this comes easily from 
$$
\sup_{|\eta-\eta^{\star}|\leq \alpha} |L''_n(\eta)-L''_n(\eta)| \leq \alpha \sup_{|\eta-\eta^{\star}|}|L_{n}^{(3)}(\eta)|
$$
where $L_{n}^{(3)}(\eta)$ is the third derivative of $L_{n}(\eta)$, and a similar handling of empirical means as before. Indeed, all functions of $\lambda$ involved are bounded as soon as $\alpha$ is such that $\eta^{\star}\geq 2 \alpha$.


\section*{Acknowledgments}

The authors would like to thank Edouard Maurel-Segala and Maxime F\'evrier for stimulating discussions
on random matrix theory and Thomas Bourgeron and Roberto Toro for having led us to study this very interesting subject
and for the discussions that we had together on genetic topics.

\bibliographystyle{abbrv}
\bibliography{biblio_anna}

\end{document}